%% file: main.tex
\documentclass[a4paper,reqno]{amsart}

\usepackage[utf8]{inputenc}
\usepackage{amsmath,amssymb,amsthm,mathtools,bm}
\usepackage[dvipsnames]{xcolor}
\usepackage[english]{babel}
\usepackage{csquotes}
\usepackage{microtype}
\usepackage[hidelinks]{hyperref}
\usepackage{fullpage}
\usepackage{tgpagella}
\usepackage{booktabs}
\usepackage[shortlabels]{enumitem}
\usepackage{tikz, pgfplots}
\usetikzlibrary{arrows.meta,calc}
\pgfplotsset{compat=1.15}
\usepackage{subcaption}

\newcommand{\R}{\mathbb{R}}
\newcommand{\N}{\mathbb{N}}
\newcommand{\cX}{\mathcal{X}}
\newcommand{\cK}{\mathcal{K}}
\newcommand{\cH}{\mathcal{H}} 
\newcommand{\Ptwo}{\mathcal{P}_2}
\newcommand{\PtwoR}{\mathcal{P}_2(\R)}
\newcommand{\Qcone}{\mathcal{Q}}
\newcommand{\Qfunc}[1]{Q\left[#1\right]}
\newcommand{\cT}{\mathcal{T}}
\newcommand{\aff}{\operatorname{aff}}
\newcommand{\conv}{\operatorname{conv}}
\DeclareMathOperator{\dist}{dist}
\DeclareMathOperator{\linspan}{span}
\newcommand{\bary}[2]{\operatorname{Bar}^{#1}({#2})}
\newcommand{\td}{\widetilde{d}}
\newcommand{\tsigma}{\widetilde{\sigma}}
\newcommand{\tLambda}{\widetilde{\Lambda}}
\newcommand{\blambda}{\bm{\lambda}}
\newcommand{\cB}{\mathcal{B}}
\newcommand{\tcB}{\widetilde{\mathcal{B}}}
\numberwithin{equation}{section}
\newcommand{\fE}{\mathfrak{E}}
\newcommand{\tfE}{\widetilde{\mathfrak{E}}}
\newcommand{\Calg}{C^{\mathrm{alg}}}
\newcommand{\Cexp}{C^{\mathrm{exp}}}
\newcommand{\alphaalg}{\alpha_{\mathrm{alg}}}
\newcommand{\alphaexp}{\alpha_{\mathrm{exp}}}
\newcommand{\stepfunc}{y}
\newcommand{\cc}{z}
\newcommand{\kk}{c}
\newcommand{\Dn}{\gamma_n}

\newtheorem{theorem}{Theorem}[section]
\newtheorem{proposition}[theorem]{Proposition}

\newtheorem{definition}[theorem]{Definition}
\newtheorem{corollary}[theorem]{Corollary}
\newtheorem{remark}[theorem]{Remark}

\newtheorem{lemma}[theorem]{Lemma}

\title{
On the convergence of 
greedy algorithms \\ in one-dimensional
Wasserstein space
}
\author{Geneviève Dusson}
\address{Université Marie et Louis Pasteur, CNRS, LmB (UMR 6623), F-25000 Besançon, France}
\email{genevieve.dusson@math.cnrs.fr}

\author{Carlo Marcati}
\address{Université Lyon 1, Centrale Lyon, INSA Lyon, Université Jean Monnet, CNRS, ICJ UMR5208, 69622 Villeurbanne, France}
\email{carlo.marcati@univ-lyon1.fr}

\author{Maximilian Penka}
\address{Université Marie et Louis Pasteur, CNRS, LmB (UMR 6623), F-25000 Besançon, France}
\email{maximilian.penka@math.cnrs.fr}

\date{\today}

\begin{document}
\begin{abstract}
We study the convergence of the greedy algorithm for approximating a
compact set by barycenters in metric spaces, in particular in the
space of one-dimensional measures endowed with the 2-Wasserstein
distance. To this end, we compare several metric $N$-widths, which measure
the best error achievable with $N$-point barycenters, using either classical
barycenters having nonnegative weights or signed barycenters, with weights of
arbitrary sign summing to one.
For the analysis we exploit that
quantile functions embed the one-dimensional 2-Wasserstein space isometrically
into a convex cone of $L^2(0,1)$, and 
barycenter sets correspond to nonlinear projections of convex or affine hulls.
We prove, both in the algebraic and exponential case,
that the convergence rate of
the greedy algorithm using signed barycenters matches the one of the metric
$N$-width defined with classical barycenters. We provide counterexamples showing
that this fails for the other combinations. 
\end{abstract}

\maketitle

\section{Introduction}
Reduced-order modeling~\cite{Quarteroni2015-db,Hesthaven2016-xm} (ROM) plays a central role for the efficient approximation of parametrized partial differential equations (PDEs).
In linear reduced-order modeling, 
the set of solutions $\mathcal{T}$ assumed to be compact
is approximated by a low-dimensional vector space, found either through a direct optimization procedure such as a proper orthogonal decomposition, or an iterative one, called a greedy algorithm.
In this context, the decay of the
best approximation error called
Kolmogorov $N$-width, 
defined for a Hilbert or Banach space $\mathcal{H}$ by
\[
    d_N(\mathcal T; \mathcal{H}) \coloneqq
    \inf_{\substack{V_N \subset \mathcal{H} \\ {\dim}(V_N) = N }}
    \adjustlimits \sup_{y\in \mathcal{T}} \inf_{v \in V_N}
    \|y - v \|_{\mathcal{H}},
\]
indicates if low-dimensional spaces can accurately approximate $\cT$.
The convergence of the greedy algorithm in Hilbert spaces
has first been studied in~\cite{buffa2012priori,binev2011convergence}, and the results have been improved and extended to Banach spaces in~\cite{DeVore2013-ec}.
The main results state that if the Kolmogorov $N$-width
decays at a certain rate (algebraic or exponential), then the greedy algorithm converges similarly, in some cases with a similar rate.

Since linear ROM performs poorly when the Kolmogorov $N$-width decays slowly, such as for transport-dominated equations~\cite{Ohlberger2015-hu},
several nonlinear reduced-order models have been recently proposed, replacing the linear vector space $V_N$ by a nonlinear manifold (see~\cite{hesthaven2026nonlinear} for a review and references therein).
Several corresponding alternative notions of widths have been proposed~\cite{cohen2023nonlinear}, such as manifold widths~\cite{devore1989optimal}, Lipschitz widths~\cite{petrova2023lipschitz}, library widths~\cite{temlyakov1998nonlinear},  $(N,M)$-widths~\cite{rim2023manifold}, and these widths are linked to the nonlinear ansatz used for approximating the solutions.

We focus in this article on metric reduced-order
modeling~\cite{Ehrlacher2020-ef,Battisti-2023,Blickhan2024Registration,
DoFeydyMula-2025,dalery2026nonlinear,Yu2026Transport}, in particular on methods where linear combinations generating
vector spaces $V_N$ above are replaced by barycenters with respect to a
well-chosen metric. In one dimension, this is strongly connected with
  ROM with cumulative distribution transforms~\cite{Ren2021RCDT,Long2025RCDT,antil2026reduced}.
In the spirit of~\cite{DeVore2013-ec}, we show how the
convergence of the greedy algorithm for metric reduced-order models relates to
the decay of the $N$-widths induced by barycenters.

Starting with the simple case of a Banach space $\cX$ 
with a norm-induced metric $\dist_{\cX}(x,y) = \|x-y\|$,
and given points $x_0,\dots,x_N \in \cX$,
a barycenter of $(x_0,\dots,x_N)$ with weights $\bm{\lambda} = (\lambda_0,\ldots,\lambda_N)\in \R_{\geq 0}^{N+1} \setminus \{0\}$ 
is defined as the minimizer of the problem
\begin{equation*}
  \bary{\bm\lambda}{x_0,\dots,x_N} := 
    \mathop{\rm argmin}_{x\in \cX} \; \sum_{i=0}^N \lambda_i\, \dist_\cX(x,x_i)^2.  
\end{equation*}
The crucial step in leaving the linear setting now is to consider a more general metric space $\cX$, but where barycenters as above are still well-defined. Note that this does not necessarily require $\cX$ to be a vector space.
In this sense, metric ROM can be seen as a generalization of linear ROM, which switches the Euclidean distance to another, well-chosen metric.

A natural notion of barycentric $N$-width inspired from the Kolmogorov $N$-width is therefore 
\begin{equation*}
    d_N(\mathcal T; \cX) \coloneqq
    \inf_{(x_i) \in \cX^{N+1}}
    \adjustlimits \sup_{y\in \mathcal{T}} \inf_{\bm{\lambda} \in \Theta_N} 
    \dist_\cX\left(y, \bary{\bm\lambda}{x_0,\dots,x_N}\right),
\end{equation*}
for which the set of weights $\Theta_{N}$ must ensure that the barycenters are well-defined. 
Candidates used in practice are either the 
$N$-dimensional simplex $\Lambda_{N}$ or the affine plane $\widetilde\Lambda_{N}$.

Again in the simple case of a norm-induced metric 
in a complete vector space, the sets of all barycenters
respectively defined with
$\Theta_{N}=\Lambda_{N}$ and $\tLambda_{N}$ 
are the convex and affine hull of $(x_0,\dots,x_N)$ 
\begin{equation*}
  \conv(x_0, \dots, x_N) = \left\{\sum_{i=0}^N\lambda_i x_i: \blambda\in \Lambda_{N}\right\}, \qquad
  \aff(x_0, \dots, x_N) = \left\{\sum_{i=0}^N\lambda_i x_i: \blambda\in \tLambda_{N}\right\}. 
\end{equation*}
Since computing best approximations involves projecting onto these sets rather than orthogonally onto a linear subspace, the results of~\cite{DeVore2013-ec} do not carry over directly and require an extended analysis.

In this paper we work in the space of probability measures with finite second moments, 
denoted $\Ptwo(\cX)$. A natural metric that lifts 
the metric of the underlying space $\cX$ is the 
2-\emph{Wasserstein distance}. 
This metric from the theory of optimal transport is 
defined as follows. Given $\mu,\nu \in \Ptwo(\cX)$, denote the set of all 
$\gamma \in \Ptwo(\cX \times \cX)$ 
that have first marginal $\mu$ and second marginal $\nu$ with $\Pi(\mu,\nu)$, and define
\begin{equation*}
    \dist_{\Ptwo}(\mu,\nu) \coloneqq \left( \inf_{\gamma \in \Pi(\mu,\nu)} \int_{\cX \times \cX} \bigl(\dist_{\cX}(x,y)\bigr)^2\ d\gamma(x,y) \right)^{1/2}.
\end{equation*}
The theory of optimal transport and metric properties of the Wasserstein distance are thoroughly discussed in many textbooks \cite{Santambrogio2015,Friesecke2024}. Finally, the metric space $(\Ptwo(\cX),\dist_{\Ptwo})$ is called Wasserstein space.

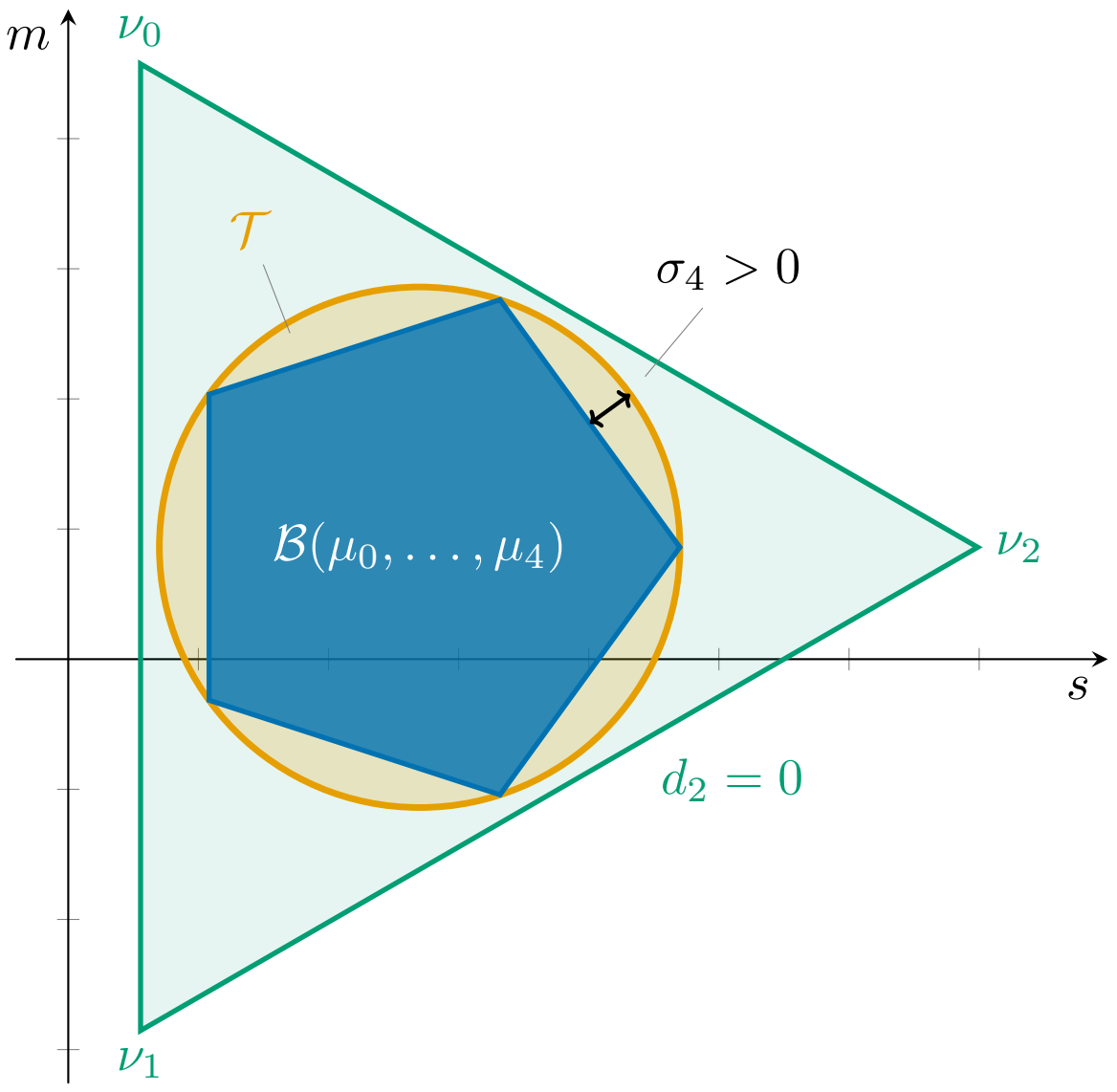
\begin{figure}[b!]
	\centering
    \input{tikz-example1}
	\caption{
    Counterexample where the manifold $\cT$ (yellow) can easily be covered
    by the convex hull of three points (green) corresponding to three Gaussians, while it is impossible to cover $\cT$ with the convex hull of any finite
    number of elements from within $\cT$ (blue).
    }
	\label{fig:example-1}
\end{figure}

\subsection{A motivating counterexample} To illustrate this, we present a motivating counterexample in the space of probability measures on $\R$ equipped with the Wasserstein metric, showing that, for the same set $\Theta_N = \Lambda_N$, the greedy error does not decay at the same rate as the metric $N$-width.
We consider a compact set of Gaussian distributions as follows:
\[
  \mathcal{T} =  \left\{ 
    \mathcal{N}(m,s^2) \,:\,
    (m - m_0)^2 + (s - s_0)^2 \le r^2
    \right\} \subset \PtwoR
\]
with $( m_0, s_0 ) \in \R\times \R_+$ and $0 < r < s_0$, visualized in yellow in $(m,s)$-coordinates in Figure~\ref{fig:example-1}.
The Wasserstein distance 
has an explicit form for Gaussian distributions
\begin{equation*}
    \dist_{\Ptwo}\left(\mathcal{N}\bigl(m,s^2\bigr), \mathcal{N}\bigl(\widetilde m,\widetilde s^2\bigr)\right)
    = \sqrt{ (m-\widetilde{m})^2 + (s-\widetilde s)^2},
\end{equation*}
i.e.~it is the Euclidean distance in the variables $(m,s).$ Therefore, the map
$\mathcal{N}(m,s^2) \mapsto (m,s)$ is an isometry into the open
$\R^2$-half-plane defined by $s > 0$.
Using that isometry, it is easy to see that $\mathcal T$ can be covered by the barycenter set with $\Theta_N = \Lambda_N$ of $\nu_0,\nu_1,\nu_2$ denoted $\mathcal{B}(\nu_0,\nu_1,\nu_2)\subset \PtwoR$ , which corresponds to the convex hull of $\nu_0,\nu_1,\nu_2$ in $(m,s)$ coordinates, and is visualized in green in Figure \ref{fig:example-1}.
However, a finite number of points of $\mathcal{T}$ are of no help in covering $\mathcal{T}$
because a convex polygon with~$N$ vertices inscribed in a disk always leaves out
a circular segment whose height is of order $N^{-2}$. 
In conclusion it holds for the $N$-width that
\begin{equation*}
    d_N(\cT; \PtwoR) = 0 \qquad\text{for all }N\ge 2,
\end{equation*}
while for all $N \in \N$ the greedy error $\sigma_N$ satisfies
\begin{equation*}
    \sigma_N\geq \adjustlimits \inf_{(\mu_i) \in \cT^{N+1}} \sup_{\nu\in \cT}
    \dist_{\Ptwo}(\nu, \cB(\mu_0,\dots,\mu_N))>\frac{c}{N^2}.
\end{equation*}
This shows the necessity 
of a careful choice for $\Theta_N$ in the considered barycenters.

\subsection{Main contributions} In this article, we thus introduce
notions of barycentric $N$-width in metric spaces and
investigate the behaviour of weak greedy algorithms for the approximation of
compact sets. 
For technical reasons,
we focus on the one-dimensional Wasserstein metric, which, using quantile functions, can be recast as a norm-induced metric (see Section~\ref{sec:2}). 
Moreover, for PDEs with transport-dominated phenomena, the Wasserstein metric is a natural choice, as translated solutions are exact barycenters between two solutions, and the corresponding metric $N$-width can be proven for $\Theta_N = \Lambda_N$ to converge faster than classical Kolmogorov $N$-width for some equations~\cite{Ehrlacher2020-ef,dalery2026nonlinear}, overcoming the so-called Kolmogorov barrier.
  We consider both $N$-widths defined with classical barycenters
  ($d_N$, Definition \ref{def:bary-n-width}, $\Theta_N = \Lambda_N$) and signed barycenters ($\td_N$,
  Definition \ref{def:bary-tilde-n-width}, $\Theta_N = \tLambda_N$) and we do the same for the errors of
  greedy sequences ($\sigma_N$ and $\tsigma_N$, Definition \ref{def:greedy}). 
  We consider whether there exist $\alpha, \beta>0$ such that the following implications are true
  \begin{align}
    \sup_{N\in\N} N^\alpha d_N <\infty \implies \sup_{N\in\N} N^\beta \sigma_N <\infty \label{eq:d<s}\\
    \sup_{N\in\N} N^\alpha d_N <\infty \implies \sup_{N\in\N} N^\beta \tsigma_N <\infty\label{eq:d<ts}\\
    \sup_{N\in\N} N^\alpha \td_N <\infty \implies \sup_{N\in\N} N^\beta \tsigma_N <\infty \label{eq:td<ts}.
  \end{align}
  We show
  \begin{enumerate}[(i)]
  \item that \eqref{eq:d<ts} is true (Theorem \ref{thm:upper}) and that the exponent $\beta=\alpha$ in the right-hand side is optimal (Proposition \ref{prop:sharp})
  \item that \eqref{eq:d<s} is false for any positive $\beta$ (Proposition \ref{prop:d<s})
  \item that \eqref{eq:td<ts} is false for any positive $\beta$ (Proposition \ref{prop:td<ts}).
  \end{enumerate}
  Note that as a consequence of \eqref{eq:td<ts} being false, we also have that
  $\sup_{N\in\N} N^\alpha\td_N <\infty$ does not imply $\sup_{N\in\N} N^\beta\sigma_N
  <\infty$ for any $\beta$, 
  since the lower bound in Proposition \ref{prop:td<ts} holds for every selection and $\cB(\mu_0,\dots,\mu_N) \subset \tcB(\mu_0,\dots,\mu_N)$.
  Furthermore, in Theorem \ref{thm:upper} we also show that if there exist $c, \alpha>0$ such that
  \begin{equation*}
     \sup_{N\in\N}\exp(cN^\alpha)d_N <\infty
  \end{equation*}
  then for all $\tilde{c}<c/(\alpha+1)$
  \begin{equation*}
    \sup_{N\in\N}\exp(\tilde{c}N^\alpha)\tsigma_N<\infty,
  \end{equation*}
  and that the exponent $\alpha$ is optimal in this case too (Proposition \ref{prop:sharp-exp}).

\subsection{Structure of the paper and notation} The rest of the article is organized as follows. In 
Section~\ref{sec:2} we properly define barycenters and present a few useful properties. In Section~\ref{sec:3} we present the definitions for the metric $N$-widths and corresponding greedy schemes.
In Section~\ref{sec:4} we present the main results comparing convergence rates, together with their proof.

In the following we use $\N_0 := \{0,1,2,\ldots \}$, and $\N := \{1,2,\ldots \}$.
The distance of a point $x \in \cX$ to a set $A \subset \cX$ is denoted as
$\dist_\cX(x,A) := \inf_{y \in A} \dist_{\cX}(x,y)$.
$P_\cX$ denotes the orthogonal projection onto $\cX$.

\section{Signed barycenters}
\label{sec:2}
In this section we provide useful definitions of barycenters for a general metric, as well as some properties in the case of the Wasserstein distance in one-dimension, using quantile functions. We start by defining signed barycenters, which correspond to the generalization of barycenters where the weights can take negative values.

To that end, we first define the two sets of weights considered hereafter
\[
  \Lambda_N \coloneqq
  \left\{   \bm\lambda \in [0,1]^{N+1} :
            \sum_{i=0}^N \lambda_i = 1 \right\},
  \qquad
  \tLambda_N \coloneqq
  \left\{   \bm\lambda \in \R^{N+1}:
            \sum_{i=0}^N \lambda_i = 1\right\}.
\]

\begin{definition}[Signed barycenter]
Let $\cX$ be a metric space, 
and let $x_0,\dots,x_N \in \cX$. A signed
barycenter in~$\cX$ 
of $(x_0,\dots,x_N)$ with weights $\bm\lambda\in \tLambda_N$ is a minimizer of the problem
\[
  \inf_{x\in \cX}\; \sum_{i=0}^N \lambda_i \dist_\cX(x,x_i)^2 .
\]
If the barycenter exists and is unique, we denote it $\bary{\bm\lambda}{x_0,\dots,x_N}$. 
The case $\bm{\lambda} \in \Lambda_N \subset \tLambda_N$ recovers the classical definition of barycenter.
\end{definition}

In all generality existence of (signed) barycenters in metric spaces is not ensured. 
But recent work establishes existence of signed barycenters in
Wasserstein space \cite{Tornabene2025}. In the following, we will work both
with $\cX = \Ptwo(\R)$, endowed with the 2-Wasserstein distance $\dist_{\Ptwo}$,
and with $\cX=L^2(0,1)$. We recall therefore two results about the existence and uniqueness of
signed barycenters that will be instrumental in the rest of the paper.

\begin{theorem}[Existence of signed barycenter,
{\cite[Th.~1 and Rem.~1]{Tornabene2025}}] \label{thm:TVS} Let $\cH$ be a separable Hilbert space and
let
$\bm\lambda \in \widetilde\Lambda_N$.
\begin{enumerate}
  \item\label{item:TVScone} Let $\cX \subset\cH$ be closed, convex, and nonempty, and let $x_0,\dots,x_N \in \cX$.
        Then $(x_i)_{i=0}^N$ admits a unique signed barycenter in $\cX$, namely
        \begin{equation*}
            \bary{\bm\lambda}{x_0,\dots,x_N}
          = P_\cX\left( {\textstyle\sum\limits_{i=0}^N \lambda_i x_i }\right),
        \end{equation*}
        where $P_\cX$ denotes the orthogonal projection onto $\cX$.
  \item Let $\mu_0,\dots,\mu_N \in \mathcal{P}_2(\cH)$. Then $(\mu_i)_{i=0}^N$ admits a signed barycenter in $(\mathcal{P}_2(\cH), \dist_{\Ptwo})$.
\end{enumerate}
\end{theorem}

\begin{remark}[Naming convention]
The existence of signed Wasserstein barycenters is established in
\cite{Tornabene2025} in a more general setting than the one we consider here;
that paper uses the term \emph{generalized} Wasserstein barycenter.
This term is overloaded by different forms of generalization in the optimal
transport literature (see, e.g., \cite{delon-2023}).
To avoid confusion, two more recent publications \cite{JZ-2026,GNT-2026} use
``signed Wasserstein barycenter'', which we also adopt in this paper.
\end{remark}

While this result guarantees existence of signed Wasserstein barycenters,
uniqueness and a characterization are hitherto unknown for general $\cH$.
The only space where it is known how to compute signed Wasserstein
barycenters is $\cH = \R$. 
The key ingredient is to identify $(\PtwoR,\dist_{\Ptwo})$ with the cone of
non-decreasing functions $\Qcone \subset L^2(0,1)$ via the \emph{quantile}
function
\begin{equation} \label{eq:quantile-mapping}
    \mu \mapsto \Qfunc{\mu}, \qquad \Qfunc{\mu}(t) := \inf\big\{x \in \R \mid t \leq
    \mu((-\infty,x]) \big\}.
\end{equation}
Vice versa, each function in $\Qcone$ defines a unique probability measure in
$\PtwoR$.
Finally, one can show~\cite{Santambrogio2015} that this bijection is an isometry and
\begin{equation*}
    \dist_{\Ptwo}(\mu,\nu) = \dist_{L^2}\bigl(\Qfunc{\mu},\Qfunc{\nu}\bigr).
\end{equation*}

\begin{corollary} \label{cor:bary-quant-character}
Let $\mu_0,\dots,\mu_N \in \PtwoR$, and let $ \Qfunc{\mu_i}$ denote their quantile
functions.
Then $\mu_0,\dots,\mu_N$ admit a unique signed barycenter with weights
$\bm\lambda \in \widetilde\Lambda_N$, and it is characterized by the quantile
function
\begin{equation*}
    \Qfunc{\bary{\bm\lambda}{\mu_0,\dots,\mu_N}} = 
    P_\Qcone\left( \sum_{i=0}^N \lambda_i \Qfunc{\mu_i} \right).
\end{equation*}
\end{corollary}

\begin{proof}
The mapping $(\PtwoR,\dist_{\Ptwo}) \to (\Qcone,\dist_{L^2}),\ \mu \mapsto \Qfunc{\mu}$ is a
bijective isometry, so it suffices to determine the signed barycenter of
$\Qfunc{\mu_0},\dots,\Qfunc{\mu_N}$ within $\Qcone$.
Since $\Qcone \subset L^2(0,1)$ is closed, convex and nonempty, Item \ref{item:TVScone} of
Theorem \ref{thm:TVS} applies with $\cH = L^2(0,1)$ and $\cX = \Qcone$.
\end{proof}

Interestingly, the projection $P_\Qcone$ can in this case be explicitly characterized~\cite[Theorem 3.1]{NatileSavare2009}. 
Indeed for $f \in L^2(0,1)$, define the integral function $F$ such that $F(t) \coloneqq
\int_0^t f(s)\,ds$ and let $G_f$ be the greatest lower semicontinuous convex
minorant of $F$. Then, the projection $P_\Qcone$ satisfies for all $t$
\begin{equation} \label{eq:isotonic-projection}
    P_\Qcone(f)(t) = \frac{d^+}{dt}G_f(t),
\end{equation}
where $d^+/dt$ indicates the right derivative.

The next two lemmas collect technical properties of the projection $P_\Qcone$.
In particular, they show that projecting a set onto $\Qcone$ does not
increase its $L^2$ distance to any element of $\Qcone$ and that $P_{\Qcone}$ maps
piecewise constant functions to piecewise constant functions, on the same partition.

\begin{lemma}
    Let $q \in \Qcone$, $A \subset L^2(0,1)$ be non-empty, and denote
    $P_{\Qcone}(A):=\{P_{\Qcone}(a) : a\in A\}$.  Then
    \begin{equation*}
        \dist_{L^2}(q,P_\Qcone(A)) \leq \dist_{L^2}(q,A).
    \end{equation*}
\end{lemma}
\begin{proof}
Since $\Qcone \subset L^2(0,1)$ is closed and convex, $P_\Qcone$ is
nonexpansive,
\[
\|P_\Qcone(f) - P_\Qcone(g)\| \le \|f-g\|
\qquad \forall f,g \in L^2(0,1) .
\]
    Moreover $P_\Qcone (q) = q$, because $q \in \Qcone$. Hence
    \[ \|q - P_\Qcone (a)\| = \|P_\Qcone (q) - P_\Qcone (a)\| \le \|q-a\| \qquad
      \forall a \in A, \]
    and taking the infimum over $a \in A$ concludes.
\end{proof}

\begin{lemma} \label{lma:projection-partition}
Let $R$ be a finite partition of the interval $(0,1)$ into subintervals and let
$V_R \subset L^2(0,1)$ be the subspace of functions that are constant on every
subinterval of $R$. Then
\begin{equation*}
    P_\Qcone(V_R) \subset V_R.
\end{equation*}
\end{lemma}
\begin{proof}
    Let $f \in V_R$. According to \eqref{eq:isotonic-projection}, in the
    continuity points of $f$,
    \[P_\Qcone(f)(t) = G'_f(t),\]
    where $G_f$ is the greatest lower semicontinuous convex minorant of $F(t) := \int_0^tf(s)\ ds$. 
    Because $f$ is constant on every subinterval of $R$, the function $F$ is
    continuous and affine on every subinterval of $R$. 
    Its greatest convex minorant is the lower boundary of the convex hull of the
    graph of $F$, and the vertices of that hull are vertices of the graph of
    $F$. Hence $G_f$ is affine on every subinterval of $R$ as well, and $G_f'$
    is constant there.
\end{proof}

 \section{\texorpdfstring{$N$}{N}-widths and greedy schemes}
 \label{sec:3}

In this section we first define $N$-widths for both classical and signed barycenters, before presenting the corresponding greedy schemes.

\subsection{Barycentric \texorpdfstring{$N$}{N}-width: classical barycenters}
The classical Kolmogorov $N$-width in a Banach space~$\cX$ measures the best
approximation error of a compact set
$\mathcal T\subset \cX$ by an $N$-dimensional linear subspace of $\cX$.
We give a definition of a compatible notion of
barycentric $N$-width in a metric space.
\begin{definition}[$N$-widths with classical barycenters]\label{def:bary-n-width}
Let $\cX$ be a metric space and let $\cT\subset \cX$ be a compact set. For $N
\in \N_0$, define
\begin{equation*}
    d_N(\mathcal T; \cX) \coloneqq
    \adjustlimits \inf_{(x_i) \in \cX^{N+1}} \sup_{y\in \mathcal{T}}
    \dist_\cX \left(y, \mathcal B(x_0,\dots,x_N)\right),
\end{equation*}
where
\[
\mathcal{B}(x_0,\dots,x_N) \coloneqq \left\{
    \bary{\bm\lambda}{x_0,\dots,x_N} : \bm\lambda \in \Lambda_{N}
\right\}
\]
is assumed to be well-defined.
\end{definition}

As mentioned in the introduction, when $\cX$ is a Banach space, this corresponds to the best approximation of
$\mathcal{T}$ by the convex hull of $N+1$ elements of $\cX$:
\begin{equation*}
  d_N(\cT; \cX) = \inf_{(x_i)\in \cX^{N+1}}\sup_{y\in\cT} \dist_\cX \left(y, \conv\{x_0, \dots, x_N\}\right).
\end{equation*}
The affine hull of $(x_i)_{i=0}^N$ contains $\conv\{x_0, \dots, x_N\}$ and it has dimension at most $N$.
The index $N$ of $d_N$ therefore counts, in this case, affine
dimensions, just as it counts linear dimensions for the classic Kolmogorov $N$-width.

In this article we are interested in the space $\PtwoR$ equipped with the 2-Wasserstein distance $\dist_{\Ptwo}$, 
in which case the barycentric $N$-width is well-defined.
Passing to quantile functions as explained in
\eqref{eq:quantile-mapping}, $d_N(\cT; \PtwoR)$ is the smallest uniform error
introduced by approximating
\begin{equation}
  \label{eq:quantile-compact}
  \cK \coloneqq \bigl\{\Qfunc{\mu}:\mu\in\cT\bigr\} \subset \Qcone
\end{equation}
with the convex hull of $N+1$ points in $\Qcone$, i.e., 
\[
    d_N(\cT; \PtwoR) = d_N(\cK; \Qcone) =
    \adjustlimits \inf_{p_0,\dots,p_N \in \Qcone} \sup_{f \in \cK}
    \dist_{L^2}\left( f, \conv\{p_0,\dots,p_N\} \right).
\]
In general, a more careful definition is needed: for
$(\mathcal{P}_2(\R^d),\dist_{\Ptwo})$ with $d > 1$ a Wasserstein barycenter might be
non-unique and for arbitrary metric spaces it might not exist at all.

\subsection{Barycentric \texorpdfstring{$N$}{N}-width: signed barycenters}
The definition of $N$-widths for signed barycenters mirrors the one that we have
given for classical barycenters with positive weights, with the only
modification that we allow for signed weights. 

\begin{definition}[$N$-widths with signed barycenters]\label{def:bary-tilde-n-width}
Let $\cX$ be a metric space and let $\cT\subset \cX$ be a compact set. For $N
\in \N_0$, define
\begin{equation*}
    \td_N(\cT; \cX) \coloneqq
    \adjustlimits \inf_{(x_i) \in \cX^{N+1}} \sup_{y\in \cT}
    \dist_\cX \bigl(y, \tcB(x_0,\dots,x_N)\bigr),
\end{equation*}
where
\[
\tcB(x_0,\dots,x_N) := \left\{
    \bary{\bm\lambda}{x_0,\dots,x_N} : \bm\lambda \in \tLambda_{N}
\right\}
\]
is assumed to be well-defined. 
\end{definition}

Note that it follows directly from Definitions~\ref{def:bary-n-width} and~\ref{def:bary-tilde-n-width} that 
$\td_N(\cT; \cX)\leq  d_N(\cT;\cX)$ 
when both $N$-widths are well defined.

When $\cX$ is a Banach space, $\td_N$ corresponds to the error of the best approximation of
$\cT$ by the affine hull of $N+1$ elements of $\cX$:
\begin{equation*}
 \td_N(\cT; \cX) = \inf_{(x_i)\in \cX^{N+1}}\sup_{y\in\cT} \dist_\cX (y, \aff\{x_0, \dots, x_N\}).
\end{equation*}
Moreover it follows from Corollary \ref{cor:bary-quant-character} that in the Wasserstein space, for a compact set $\cT\subset\PtwoR$, one can use quantile functions and obtain
\[
    \td_N(\cT; \PtwoR) = 
    \adjustlimits \inf_{p_0,\dots,p_N \in \Qcone} \sup_{f \in \cK}
    \dist_{L^2}\left( f, P_\Qcone(\aff\{p_0,\dots,p_N\}) \right),
\]
with $\cK$ defined as in \eqref{eq:quantile-compact}.

\subsection{Greedy schemes}
We now define the greedy algorithm whose convergence rates we compare with barycentric $N$-widths hereafter.
\begin{definition}[Greedy sequence]
 \label{def:greedy}
Let $\cT \subset \PtwoR$ be compact and nonempty. Let $\Theta_N$ denote either
$\Lambda_N$ or $\tLambda_N$.
Define, for $N\in\N_0$ and $(\mu_i)_{i=0}^{N}\subset \cT$, the error functional
\begin{equation*}
     E_N(\nu)\coloneqq \inf_{\blambda \in \Theta_{N}}
    \dist_{\Ptwo}\left(\nu, \bary{\blambda}{\mu_0, \dots, \mu_{N}} \right),
\end{equation*}
Fix $\gamma\in(0,1]$. 
A sequence $( \mu_i )_{i\in\N_0} \subset \cT$ is a
weak greedy sequence for $\cT$ (with constant $\gamma$) if for all $N \in \N_0$
\begin{equation*}
  E_N(\mu_{N+1}) \geq \gamma \max_{\nu\in\cT} E_N(\nu).
\end{equation*}
We denote the error at the $N$th step by
\begin{equation*}
  \max_{\nu\in\cT}E_N(\nu) \eqqcolon
  \begin{cases}
    \sigma_N&\text{if }\Theta_{N} = \Lambda_{N}\\
    \tsigma_N&\text{if }\Theta_{N} = \tLambda_{N}.
  \end{cases}
\end{equation*}
We refer to the sequences as weak greedy sequences with classical or signed
barycenters, respectively.
\end{definition}
It follows directly from this definition that the sequences $\sigma_N$ and $\tsigma_N$ are nonincreasing in $N$.

For the Wasserstein space $\PtwoR$, 
passing to quantile functions, we obtain an equivalent characterization of a
weak greedy sequence through quantile functions. Denote $q_i:=\Qfunc{\mu_i}$ for
$i\in \N_0$. Then, let
\begin{align}
\fE_N(f) := \dist_{L^2} \Big( f, \conv\{q_0,\ldots,q_{N}\}  \Big), 
\\
\tfE_N(f) := \dist_{L^2} \Big( f, P_{\Qcone} \big( \aff\{q_0,\ldots,q_{N}\} \big) \Big).
\end{align}
For every $N \in \N_0$, the \emph{weak greedy rule} with, respectively, classical
and signed barycenters, selects $q_{N+1} \in \cK:= \{\Qfunc{\mu}
: \mu \in \cT\}$
such that
\begin{subequations}
  \label{eq:greedy-rule}
\begin{align}
    \fE_N(q_{N+1}) \ge \gamma\max_{f\in\cK}\fE_N(f), 
    \\
    \tfE_N(q_{N+1}) \ge \gamma\max_{f\in\cK}\tfE_N(f). 
\end{align}
\end{subequations}
It also follows that $\sigma_N = \max_{f\in\cK}\fE_N(f)$ and $\tsigma_N = \max_{f\in\cK}\tfE_N(f)$.

\section{Comparison of rates of convergence}
  \label{sec:4}
In this section, we provide the main results of this article. We start by showing the
positive result, i.e., that implication \eqref{eq:d<ts} is true. Then, in
Sections \ref{sec:d<s} and \ref{sec:td<ts}, we construct counterexamples to
implications~\eqref{eq:d<s} and \eqref{eq:td<ts}, respectively.
\subsection{Comparison of $d_N$ and $\tsigma_N$}
We start by comparing the $N$-width defined with classical barycenters and the
error of the greedy scheme that uses signed barycenters.

\begin{theorem}[Algebraic greedy rate]\label{thm:upper}
Let $\cT\subset \PtwoR$ be a compact set, where $\PtwoR$ is equipped with the
Wasserstein distance.
Let $\tsigma_N$ be the error of a weak greedy sequence for $\cT$ defined as in Definition~\ref{def:greedy} with $\Theta_N=\tLambda_N$.
If, for some $\Calg_1>0$ and $\alphaalg>0$, 
\begin{equation} \label{eq:width-rate}
    d_N(\cT; \PtwoR)\leq \Calg_1 (N+1)^{-\alphaalg}, \qquad \forall N \in \N_0 ,
\end{equation}
then
\begin{equation} \label{eq:greedy-rule-rate}
    \tsigma_N \leq \Calg_2 (N+1)^{-\alphaalg},  \qquad \forall N \in \N_0.
\end{equation}
If, for some $\Cexp_1, c_1, \alphaexp>0$, 
\begin{equation} \label{eq:width-rate-exp}
    d_N(\cT; \PtwoR)\leq \Cexp_1\exp(-c_1 (N+1)^{\alphaexp}), \qquad \forall N \in \N_0,
\end{equation}
then
\begin{equation} \label{eq:greedy-rule-rate-exp}
    \tsigma_N\leq \Cexp_2 \exp(-c_2 (N+1)^{\alphaexp}),  \qquad \forall N \in \N_0.
\end{equation}
Admissible choices for the constants are
\begin{equation} \label{eq:C2}
    \Calg_2 = \frac{2^{\alphaalg + 3/2}e^{\alphaalg+1/2}}{\gamma}\Calg_1, 
    \qquad 0 < c_2 < \frac{c_1}{\alphaexp+1}, \qquad 
    \Cexp_2 = \frac{2e^{1/2}}{\gamma\sqrt{1-((\alphaexp+1)c_2/c_1)^{1/\alphaexp}}}\Cexp_1.
\end{equation}
\end{theorem}

\subsubsection{Proof of Theorem \ref{thm:upper}}
For the proof we start by collecting two lemmas that hold for an arbitrary Hilbert
space $\cH$.
In both Lemmas, we consider an operator $T: \R^N \to \cH$ whose ``columns'' are
the elements $x_1, \dots, x_N \in \cH$, so that for the Hilbert--Schmidt norm
\begin{equation*}
\|T\|_{\mathrm{HS}}^2 = \sum_{j=1}^N\|x_j\|^2_\cH.
\end{equation*}

The first lemma links the distance between one ``column'' of $T$ and the
subspace generated by the previous ones to the singular values of $T$.

\begin{lemma}\label{lem:gram}
Let $x_1,\ldots,x_N\in \cH$, and let $T:\R^N\to \cH$ be defined by $Te_j=x_j$,
where $(e_j)$ denotes the standard basis of $\R^N$.
Set $r_1 = \|x_1\|$ and for $j = 2,\dots, N$
\[
 r_j = \dist_\cH \bigl( x_j,
       \linspan\{x_1,\ldots,x_{j-1}\}
       \bigr) .
\]
Let $s_1,\ldots,s_N$ denote the singular values of $T$. Then
\begin{equation}
 \prod_{j=1}^N r_j^2=\prod_{j=1}^N s_j^2.
\end{equation}
\end{lemma}

\begin{proof}
If the vectors $x_j$ are linearly dependent, both products vanish.  
Otherwise, Gram--Schmidt orthogonalization yields $T=UR$,
where $R$ is an upper triangular $N\times N$-matrix with diagonal entries $r_1,\ldots,r_N$.  
Therefore,
\[
 \prod_{j=1}^N s_j^2 = \det(T^*T)=\det(R^*R)=\prod_{j=1}^N r_j^2.
\]
\end{proof}

In the second lemma we estimate the individual singular values of the operator $T$.

\begin{lemma}\label{lem:multiscale}
Let $x_1,\ldots,x_N\in \cH$, and define $T:\R^N\to \cH$ via $Te_j=x_j$, with singular values
arranged in nonincreasing order $s_1\geq \dots \geq s_N$.
Let, for $k=0,\ldots,N-1$, $W_{k} \subset \cH$ be a subspace of dimension at most $k$.
Then, for all $j\in\{1, \dots,N\}$, and $\beta\in(0,1)$ the singular values of $T$ satisfy
\begin{equation}
    s^2_j\ \leq\ \frac{N}{(1-\beta)j} \, 
    \max_{i=1, \dots, N}\dist_{\cH}\left(x_i, W_{\lfloor  \beta j \rfloor}\right)^2.
\end{equation}
\end{lemma}

\begin{proof}
Fix $k \in \{0,\dots,N-1\}$ and note that $P_{W_{k}}T$ has rank at most $k$. 
Standard low-rank approximation theory, followed by the definition of $T$, give
\begin{equation}
 \label{eq:tail-algebraic}
 \sum_{i=k+1}^N s_i^2\leq \|T-P_{W_{k}}T\|_{\mathrm{HS}}^2 
 =\sum_{i=1}^N \|x_i-P_{W_{k}}x_i\|_\cH^2
 =\sum_{i=1}^N\dist_\cH(x_i,W_{k})^2.
\end{equation}
Fix $j\in\{1,\ldots,N\}$ and choose
$k=\lfloor \beta j \rfloor$.
Since the singular values are nonincreasing,
\[
 (j-k)s_j^2
 \leq\sum_{i=k+1}^j s_i^2
 \leq\sum_{i=k+1}^N s_i^2,
\]
and the conclusion follows from~\eqref{eq:tail-algebraic} and $j-k \geq (1-\beta)j$. 
\end{proof}

\begin{proof}[Proof of Theorem \ref{thm:upper}]
Let $q_0,q_1,q_2,\ldots \in \cK = \{\Qfunc{\mu}:\mu\in \cT\}$ be a sequence of 
quantile functions of measures selected by the weak
greedy rule \eqref{eq:greedy-rule}.
Define the difference vectors $x_j := q_j - q_0$ for $j \in \N$
and consider the sequence of subspaces spanned by them:
\begin{equation*}
    V_{j} := \linspan\{x_1,\dots,x_j\}
    = \left\{ q - q_0 : q \in \aff\{q_0,\dots,q_{j}\} \right\},
\end{equation*}
and $V_0 := \{0\}$.
Then, using the nonexpansiveness of $P_\Qcone$,
\begin{equation}\label{eq:xj-Vj}\begin{aligned}
    \dist_{L^2}(x_{j+1},V_{j})
    &= \dist_{L^2}(q_{j+1},\aff\{q_0,\dots,q_{j}\})\\
    &\ge \dist_{L^2}(q_{j+1},P_{\Qcone}(\aff\{q_0,\dots,q_{j}\})) \\
    &= \tfE_{j}(q_{j+1}) \ge \gamma \tsigma_{j}.
\end{aligned}\end{equation}
For any $N \in \N_0$ let $T:\R^{N+1} \to L^2(0,1)$ be defined by $Te_j = x_j$ for all $j = 1,\dots,N+1$, with singular values $s_1 \geq \dots \geq s_{N+1}$. Using the monotonicity of $\tsigma_j$,~\eqref{eq:xj-Vj}, 
and Lemma \ref{lem:gram},
we obtain
\begin{equation} \label{eq:use-lem-gram}
     (\gamma \tsigma_{N})^{2(N+1)}
     \le \prod_{j=1}^{N+1} \dist(x_{j},V_{j-1})^{2}
     = \prod_{j=1}^{N+1} s_j^{2}.
\end{equation}

In order to apply Lemma \ref{lem:multiscale} to the singular values on the right hand side,
we need to construct the spaces $W_k$.
For all $k\in\N_0$ consider the linear space $W_k$ realizing 
\cite[Th.~2.3]{Pinkus1985} the minimum
\begin{equation}
\label{eq:Wk}
    \min_{\substack{W\subset L^2\\\dim(W)\leq k}} \sup_{i \in \N} \dist_{L^2}(q_i-q_0, W) =\sup_{i \in \N} \dist_{L^2}(q_i-q_0, W_k).
\end{equation}
Next, we need to link the distances on the right hand side to the $N$-width. 
Fix any $p_0, \dots, p_k \in \Qcone$ and denote by $z(q)$ the projection of $q\in\Qcone$ onto
$\conv\{p_0, \dots, p_k\}$. Since $z(q_i)-z(q_0)\in \linspan\{p_1-p_0, \dots, p_k-p_0\}$, we have that, for all $i \in \N$,
\begin{equation}
  \label{eq:spandiff}
  \begin{multlined}
    \dist_{L^2}(q_i-q_0,\linspan\{p_1-p_0,\dots,p_k-p_0\}) 
    \leq \|q_i - q_0 - (z(q_i) - z(q_0))\|
    \\
    \leq \|q_i-z(q_i)\| + \|q_0 - z(q_0)\| 
    \leq 2\sup_{q \in \cK} \dist(q,\conv\{p_0,\dots,p_k\}).
  \end{multlined}
\end{equation}
Because the span of $k$ elements is at most $k$-dimensional, using~\eqref{eq:spandiff} and the definition of $W_k$~\eqref{eq:Wk},
\begin{multline*}
    \sup_{i \in \N} \dist_{L^2}(\underbrace{q_i - q_0}_{x_i},W_k) 
    = \min_{\substack{W \subset L^2 \\ \dim W \leq k}} \sup_{i \in \N} \dist_{L^2}(q_i - q_0,W) 
    \\
    \leq \sup_{i \in \N} \dist_{L^2}(q_i-q_0,\linspan\{p_1-p_0,\dots,p_k-p_0\})
    \leq 2\sup_{q\in\cK} \dist_{L^2}(q, \conv\{p_0, \dots, p_k\}).
\end{multline*}
The left-hand side does not depend on $p_0, \dots, p_k$, hence we can take
the infimum of the right-hand side over all $(p_0, \dots, p_k) \in \Qcone^{k+1}$,
to obtain
\begin{equation*}
    \sup_{i\in \N} \dist_{L^2}(x_i, W_k)
    \leq 2 \adjustlimits \inf_{p_0,\dots,p_k \in \Qcone} \sup_{q \in\cK} \dist_{L^2}(q,\conv\{p_0,\dots,p_k\})
    = 2 d_k (\cT;\PtwoR).
\end{equation*}
Using \eqref{eq:use-lem-gram} and Lemma \ref{lem:multiscale} we obtain for any $\beta \in (0,1)$ 
\begin{equation} \label{eq:use-lem-multiscale}
    (\gamma \tsigma_{N})^{2(N+1)}
    \leq \prod_{j=1}^{N+1} s_j^{2}
    \leq \left(\frac{N+1}{1-\beta}\right)^{N+1} \prod_{j=1}^{N+1} \frac{1}{j}\max_{i=1, \dots, N+1}\dist_{L^2}(x_i, W_{\lfloor \beta j \rfloor})^2.
\end{equation}
Under assumption \eqref{eq:width-rate}, choosing $\beta = 1/2$ so that $k = \lfloor j/2 \rfloor$,
\[
d_{\lfloor j/2 \rfloor}(\cT;\PtwoR) \leq \Calg_1 (j/2)^{-\alphaalg},
\]
and using Stirling's inequality $N! \geq (N/e)^N$,
\begin{align*}
  \prod_{j=1}^{N+1} \frac{1}{j}\max_{i=1, \dots, N+1}\dist_{L^2}\left(x_i, W_{\lfloor j/2 \rfloor}\right)^{2}
  \leq \left( 2^{\alphaalg+1}  \Calg_1\right)^{2(N+1)}  
  \left(\frac{e}{N+1}\right)^{(N+1)(2\alphaalg+1)}.
\end{align*}
Hence, plugging this with $\beta = 1/2$ in \eqref{eq:use-lem-multiscale}, 
and taking the $2\,(N+1)$-th root gives \eqref{eq:greedy-rule-rate}.

Under assumption \eqref{eq:width-rate-exp}, for any $\beta \in (0,1)$ 
\[
d_{\lfloor \beta j \rfloor}(\cT; \PtwoR) \leq \Cexp_1 \exp(-c_1 (\beta j)^{\alphaexp}),
\]
and therefore
\begin{align*}
  \prod_{j=1}^{N+1} \frac{1}{j}\max_{i=1, \dots, N+1} \dist_{L^2}\left(x_i, W_{\lfloor \beta j \rfloor}\right)^{2}
  & \leq
  4^{N+1}(\Cexp_1)^{2(N+1)}  
  \prod_{j=1}^{N+1} \frac{1}{j}\exp(-2c_1(\beta j)^{\alphaexp})
  \\ & \leq
  \left(2\,\Cexp_1\right)^{2(N+1)}
  \left(\frac{e}{N+1}\right)^{N+1} \exp\left(-2c_1\beta^{\alphaexp}\frac{(N+1)^{\alphaexp+1}}{\alphaexp+1}\right).
\end{align*}
Plugged into \eqref{eq:use-lem-multiscale}, this yields
\begin{equation*}
    \tsigma_{N} \leq \frac{2e^{1/2}}{\gamma\sqrt{1-\beta}}\Cexp_1\exp(-c_1 \beta^{\alphaexp} (N+1)^{\alphaexp}/(\alphaexp+1)), 
\end{equation*}
and setting $c_2 = c_1 \beta^{\alphaexp}/(\alphaexp+1)$, i.e. $\beta=((\alphaexp+1)\,c_2/c_1)^{1/\alphaexp}$, 
gives the last relation in \eqref{eq:C2}.
\end{proof}

\begin{remark}
This proof of Theorem \ref{thm:upper} has certain similarity with the techniques of 
\cite{DeVore2013-ec} in encoding the greedy sequence in an operator $T$. Compared to 
that reference, we choose a slightly more direct way in linking the greedy error to the 
singular values of $T$ and account for working in affine coordinates and the nonlinear 
projection $P_\Qcone$.

Vice versa, our proof never utilizes any specific features of $L^2(0,1)$ besides being a
Hilbert space and $\Qcone$ being closed and convex, so that all arguments apply for the 
choice $\Qcone = \cH$ as well. Adapting for the linear span in the Kolmogorov $N$-width and 
greedy error leads to the known relations of \cite{DeVore2013-ec} but with slightly different constants.
\end{remark}

\subsubsection{Sharpness of the convergence rates}
  \label{sec:sharpness}
In this section we show that the exponents $\alphaalg$ and $\alphaexp$ in Theorem~\ref{thm:upper} cannot be improved. To do so,
we will consider quantile functions and construct a compact set for which the
exponent is indeed sharp.
Let $(c_n)_{n\in\N}\subset \R_+$ such that $c_n \to 0$.
Let $t_n \coloneqq 1 - 2^{-n}$ 
and consider the family of step functions
\begin{equation} \label{eq:sharp-family}
    \stepfunc_0 \coloneqq 0, \qquad
    \stepfunc_n \coloneqq c_n2^{n/2}\,\bm{1}_{[t_n,1]}\, , \quad n \ge 1.
\end{equation}
Those are the quantile functions of the measures
\begin{equation} \label{eq:sharp-family-measures}
    \mu_0 = \delta_0, \qquad \mu_n = t_n\delta_0 +
    (1-t_n)\delta_{c_n 2^{n/2}}\, , \quad n \ge 1.
\end{equation}

The properties of this family used in the subsequent arguments are summarized in the following lemma, which is proved in Appendix~\ref{sec:step-fun}.

\begin{lemma} \label{lem:spike-fam-measure}
Given a positive sequence $(c_n)$ with $c_n \to 0$, let $\cT$ denote the family of measures \eqref{eq:sharp-family-measures}.
\begin{enumerate} \renewcommand{\labelenumi}{(\roman{enumi})}
    \item The family $\cT \subset \PtwoR$ is compact in the Wasserstein metric.
    \item For all $N \in \N_0$
    \begin{equation} \label{eq:spike-width-measure}
        d_N(\cT; \PtwoR) \leq \max_{n \geq N+1} c_n.
    \end{equation}
    \item For $N \geq 0$ and every selection $\nu_0,\dots,\nu_N \in \cT$,
    \begin{equation} \label{eq:spike-lower-measure}
        \max_{n \geq 0} \dist_{\Ptwo}\left(\mu_n, \tcB(\nu_0,\dots,\nu_N)\right) \geq 3^{-1/2}\min_{n=1,\dots,N+2} c_n .
    \end{equation}
\end{enumerate}
\end{lemma}

The following proposition exhibits that the optimal universally guaranteed
algebraic exponent from Theorem \ref{thm:upper} is precisely $\alpha$.
Recall that the greedy error $\tsigma_N$ is measured against the $N+1$ selected
elements $\mu_0,\dots,\mu_N$. Accordingly, the lower bound below is stated for
an arbitrary selection of $N+1$ elements.

\begin{proposition}
  \label{prop:sharp}
  Let $\alpha > 0$. Let
  $\cT_\alpha :=\{\mu_n: n \in \N_0\}$
  be the family \eqref{eq:sharp-family-measures} 
  with $c_n = n^{-\alpha}$. Then the set $\cT_\alpha \subset \PtwoR$ 
  is compact and satisfies
  \begin{equation} \label{eq:sharp-width}
    d_N(\cT_\alpha; \PtwoR) \le (N+1)^{-\alpha}, \qquad N \ge 0.
  \end{equation}
  Moreover, for every $N \ge 0$ and every selection $\nu_0,\dots,\nu_{N} \in
  \cT_\alpha$,
  \begin{equation} \label{eq:sharp-lower}
    \max_{\mu \in \cT_\alpha}
    \dist_{\Ptwo}(\mu,\tcB(\nu_0, \dots, \nu_N))
    \,\geq\, 3^{-1/2} (N+2)^{-\alpha}.
  \end{equation}
  In particular, for every $\gamma \in (0,1]$, every weak greedy sequence for
  $\cT_\alpha$ selected according to Definition \ref{def:greedy} satisfies
  \begin{equation} \label{eq:sharp-sigma}
    \frac{2^{-\alpha}}{\sqrt{3}}(N+1)^{-\alpha} \le \tsigma_N \le 
    \frac{2^{\alpha+3/2}e^{\alpha+1/2}}{\gamma}(N+1)^{-\alpha}, \qquad N \ge 0.
  \end{equation}
  Hence the exponent in $(N+1)^{-\alpha}$ in \eqref{eq:greedy-rule-rate} cannot be
  universally improved.
\end{proposition}

\begin{proof}
According to Lemma \ref{lem:spike-fam-measure}, $\cT_\alpha$ is compact. Moreover since $c_n$ is decreasing,
\eqref{eq:sharp-width} is \eqref{eq:spike-width-measure} and 
\eqref{eq:sharp-lower} is \eqref{eq:spike-lower-measure}.
Finally, let $\nu_0,\dots,\nu_N$ be a weak greedy sequence for
$\cT_\alpha$.
Since $N+2 \leq 2(N+1)$ for $N \geq 0$, the lower bound in \eqref{eq:sharp-sigma}
follows from \eqref{eq:sharp-lower},
while the upper bound follows from \eqref{eq:sharp-width} and
Theorem \ref{thm:upper} applied with $\Calg_1 = 1$.
\end{proof}

We now turn to prove the exponent sharpness in the exponential case.

\begin{proposition} \label{prop:sharp-exp}
  Let $\alpha, c > 0$. Let
  $\cT_\alpha =\{\mu_n: n \in \N_0\}$
  be the family \eqref{eq:sharp-family-measures}
  with $c_n = \exp(-c n^\alpha)$. 
  Then the set $\cT_\alpha$ is compact and satisfies
  \begin{equation}
    d_N(\cT_\alpha;\PtwoR) \le \exp(-c\,(N+1)^{\alpha}), \qquad N \ge 0.
  \end{equation}
  Moreover, for every $N \geq 0$ and every selection $\nu_0,\dots,\nu_N \in
  \cT_\alpha$,
  \begin{equation}
    \max_{\mu \in \cT_\alpha}
    \dist_{\Ptwo}(\mu,\tcB(\nu_0, \dots, \nu_N))
    \,\geq\, 3^{-1/2} \exp(-c\,(N+2)^{\alpha}).
  \end{equation}
  In particular, for every $\gamma \in (0,1]$, every weak greedy sequence for
  $\cT_\alpha$ selected according to Definition \ref{def:greedy} satisfies
  \begin{equation}
    \frac{1}{\sqrt{3}}\exp(-2^\alpha c\,(N+1)^{\alpha}) 
    \le \tsigma_N \le 
    \frac{2^{3/2}e^{1/2}}{\gamma}
    \exp\left(\frac{-c\,(N+1)^{\alpha}}{2^\alpha(\alpha+1)}\right), \qquad N \ge 0.
  \end{equation}
  Hence the exponent $\alpha$ in $\exp\left(-c\,(N+1)^{\alpha}\right)$ in 
  \eqref{eq:greedy-rule-rate-exp} cannot be universally improved.
\end{proposition}
\begin{proof}
  The proof is the same as that of Proposition \ref{prop:sharp}.
  Here we have made a specific choice of $c_2$ when applying
  \eqref{eq:greedy-rule-rate-exp}, such that
  $((\alphaexp+1)c_2/c_1)^{1/\alphaexp} = 1/2$.
\end{proof}

\subsection{Comparison of $d_N$ and $\sigma_N$}
  \label{sec:d<s}
The example of the Gaussian family from the introduction already showed that $d_N$ and $\sigma_N$ might not share an algebraic relation. In the following, we work out a more dedicated counter-example showing that \eqref{eq:d<s} is false for any positive $\beta$, where the Gaussian example could only provide the case $\beta > 2$.

We need a first technical lemma before providing the counter example.
\begin{lemma} 
  \label{lemma:OrderedBall}
  Let $n \ge 2$ and let $\cc \in \R^n$ be defined by $\cc^{(j)} = \sqrt{2}(j-1)$ for $j=1,\dots,n$. 
  Consider the cone of vectors with ordered components
  \begin{equation*}
    A = \left\{a\in \R^n: a^{(1)} \leq a^{(2)}\leq \dots \leq a^{(n)}\right\}.
  \end{equation*}
  Then, the closed Euclidean unit ball centered at $\cc$ satisfies $\overline{B}(\cc)\subset A$.
  Furthermore, there exist $n+1$ points $p_1, \dots, p_{n+1}\in A$ such that
  $\overline{B}(\cc)\subset \conv\{p_1, \dots, p_{n+1}\}$. Finally, there exists
  a constant $\kk > 0$ that depends only on $n$ such that for all $N \in \N$ and
  all $\{b_1,\dots,b_N\} \subset \overline{B}(\cc)$, it holds that
  \begin{equation} \label{eq:lowerbound-ball}
  \max_{x\in \overline{B}(\cc)}\dist_{\ell^2}(x, \conv\{b_1, \dots, b_N\})\geq \kk N^{-2/(n-1)}.
  \end{equation}
\end{lemma}
\begin{proof}
  We start by showing that $\overline{B}(\cc)\subset A$: For all $x\in
  \overline{B}(\cc)$, one has 
  $\left(x^{(j+1)} -\cc^{(j+1)}\right)^2 +  \left(x^{(j)}-\cc^{(j)}\right)^2\leq 1$,
  hence 
  $\left|x^{(j+1)} - \cc^{(j+1)}\right| + \left|x^{(j)} -\cc^{(j)}\right|\leq \sqrt{2}$ 
  by Cauchy--Schwarz and
  $x^{(j+1)} - x^{(j)} \geq \cc^{(j+1)}-\cc^{(j)} - \sqrt{2} = 0$.
  Then, we remark that $\overline{B}(\cc)$ is contained in
  \begin{equation*}
    \widetilde{A}\coloneqq\{a\in \R^n: m\coloneqq -1 \leq a^{(1)} \leq \dots \leq a^{(n)} \leq \sqrt{2}(n-1)+ 1\eqqcolon M\} \subset A,
  \end{equation*}
   which is the convex hull of the points 
  \begin{equation*}
    p_j = (\underbrace{m , \dots ,m}_{j-1 \text{ times}}, \underbrace{M, \dots, M}_{n-j +1\text{ times}}),\qquad j=1, \dots, n+1.
  \end{equation*}
  Finally, the last relationship \eqref{eq:lowerbound-ball} is the approximation of a ball by an inscribed
  simplex \cite{BaranyFuredi1988}.
\end{proof}

\begin{proposition}\label{prop:d<s}
 For all $\beta>0$, there exist a compact set $\cT\subset \PtwoR$ and a constant
 $\kk>0$ such that $d_N(\cT; \PtwoR) = 0$ for all integers
 $N\geq 2/\beta+1$ and
for every weak greedy sequence with classical barycenters
 \begin{equation}
   \label{eq:lowerbound-d<s}
    \sigma_N \geq \adjustlimits \inf_{(\mu_i) \in \cT^{N+1}} \sup_{\nu\in \cT}
    \dist_{\Ptwo}(\nu, \cB(\mu_0,\dots,\mu_{N})) \geq \kk\,(N+1)^{-\beta} \qquad\text{for all }N\in\N_0.
 \end{equation}
\end{proposition}
\begin{proof}
  Let $n\coloneqq\lceil 2/\beta\rceil+1$, so that $n \ge 2$.
  Let $A, \cc, p_1,\dots,p_{n+1}$ and $\kk > 0$ be as in Lemma \ref{lemma:OrderedBall}, and for $a \in A$ consider the measure
  \begin{equation*}
      \mu_a \coloneqq \frac{1}{n}\sum_{j=1}^n \delta_{\sqrt{n}a^{(j)}} \in \PtwoR,
  \end{equation*}
  with corresponding quantile function given by $\Qfunc{\mu_a}(t) = \sqrt{n}a^{(j)}$ for $t \in \left(\frac{j-1}{n},\frac{j}{n}\right), j = 1,\dots,n$.
  In particular, $a \mapsto \Qfunc{\mu_a}$ is linear on $A$ and
  \begin{equation}
    \label{eq:isometryOrdered}
    \dist_{\Ptwo}(\mu_a, \mu_b)^2 = \| \Qfunc{\mu_a} - \Qfunc{\mu_b}\|^2_{L^2}
    = |a-b|^2 \qquad \text{for all } a,b \in A.
  \end{equation}
  Define $\cT = \left\{ \mu_a : a\in \overline{B}(\cc)\right\}$, which is compact by
  \eqref{eq:isometryOrdered}.
  Since $\overline{B}(\cc) \subset \conv\{p_1,\dots,p_{n+1}\}$, we identify $\cT \subset \cB(\mu_{p_1}, \dots, \mu_{p_{n+1}})$, and obtain $d_n(\cT;\PtwoR) =0$.
  Since $d_N$ is non-increasing in $N$, we also obtain $d_N(\cT;\PtwoR) = 0$ for all $N \geq n$.
  
  Finally, let $N \in \N_0$ and let $\nu_0,\dots,\nu_{N}\in\cT$ be arbitrary, 
  with $\nu_{i-1} = \mu_{b_i}$ for some $b_1,\dots,b_{N+1} \in \overline{B}(\cc)$.
  By \eqref{eq:isometryOrdered}, for every $y \in \overline{B}(\cc)$,
	\begin{equation*}
		\dist_{\Ptwo}(\mu_y,\cB(\nu_{0}, \dots , \nu_{N})) = \dist_{\ell^2}(y,\conv\{b_1,\dots,b_{N+1}\}).
	\end{equation*}
	Taking the supremum over $\nu = \mu_y\in \cT$ and using \eqref{eq:lowerbound-ball}, 
    the left-hand side is at least $\kk(N+1)^{-2/(n-1)} \geq \kk(N+1)^{-\beta}$, by the choice of $n$. 
	As the points were arbitrary, taking the infimum yields the second
	inequality in \eqref{eq:lowerbound-d<s}. 
    The first one holds because the first $N+1$ elements of the greedy sequence lie in~$\cT$.
\end{proof}

\subsection{Comparison of $\td_N$ and $\tsigma_N$}
  \label{sec:td<ts}
Since the previous example shows that it does not suffice to consider ordinary
barycenters in the greedy approximation, one could ask in turn whether it would
not be more natural to consider signed barycenters in the definition of the
barycenter $N$-width as well; that is, recalling Definition \ref{def:bary-tilde-n-width},
\[
    \td_N(\cT;\PtwoR) = \td_N(\cK;\Qcone) =
    \adjustlimits\inf_{p_0,\dots,p_N \in \Qcone} \sup_{f \in \cK}
    \dist_{L^2}(f,P_\Qcone (\aff\{p_0,\dots,p_N\})).
\]
In the following proposition, we construct a family of quantile 
functions and show that the nonlinear projection $P_\Qcone$ can 
generate the whole family from the affine hull of a finite set of 
generators outside the family.
This family of quantile functions cannot instead be generated by 
applying $P_\Qcone$ to affine combinations of a finite selection of 
generators within the family.
The set is built from studying the same family of discrete measures 
introduced in Section \ref{sec:sharpness}.
The relevant feature of the set constructed in the proposition
is that the breakpoints $( t_n )_{n\in \N_0}$ of the quantile
functions are countable and distinct. Any finite, affine combination 
of quantile functions from this set has a finite number of 
breakpoints, and Lemma \ref{lma:projection-partition} prevents the 
projection from creating any additional one, thus giving a lower bound 
on $\tsigma_N$.

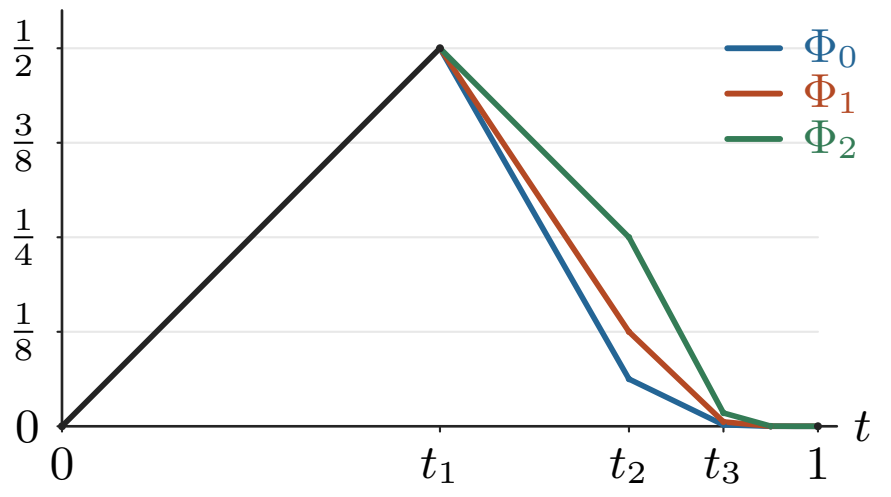
\begin{figure}[b]
    \centering
    \vspace{0pt}
    \tikzset{every picture/.append style={scale=0.50}}
    \input{phi-functions.tex}
    \caption{Functions $\Phi_0,\Phi_1,\Phi_2$ defined in~\eqref{eq:Phi-knots}}
    \label{fig:phi}
\end{figure}

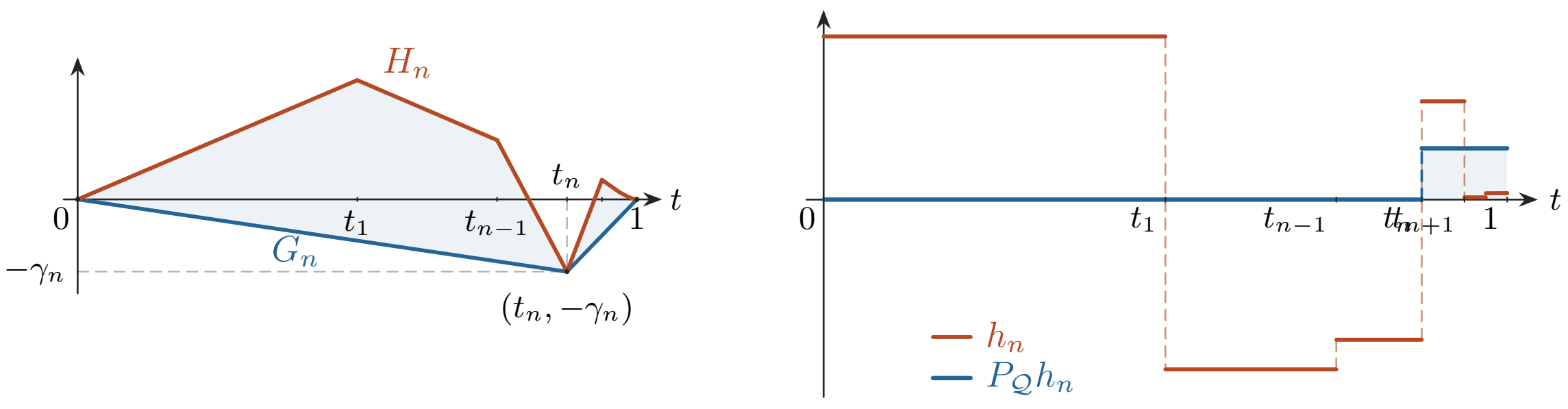
\begin{figure}[t]
  \centering
  \begin{subfigure}[t]{0.49\textwidth}
    \centering
    \tikzset{every picture/.append style={scale=0.45}}
    \raisebox{19pt}{\input{hn-gn.tex}}
    \caption{$H_n$ and $G_n$, $y$-axis not to scale.}\label{subfig:HnGn}
  \end{subfigure} 
  \hfill
  \begin{subfigure}[t]{0.49\textwidth}
    \centering
    \tikzset{every picture/.append style={scale=0.55}}
    \input{hn-projection.tex}
    \caption{$h_n$ and $P_{\Qcone}h_n$, $y$-axis not to scale.}\label{subfig:hn-Phn}
  \end{subfigure}
  \caption{Schematic representation of the functions constructed in the proof of
  Proposition \ref{prop:td<ts}. }
\end{figure}

\begin{proposition}\label{prop:td<ts}
There exists a compact set $\cT\subset \PtwoR$ such that $\td_N(\cT; \PtwoR) = 0$ for 
$N\geq 5$ and, for every weak greedy sequence with signed barycenters,
 \begin{equation}
 \label{eq:sigmaNlowerbound}
    \tsigma_N\geq \adjustlimits \inf_{(\mu_i) \in \cT^{N+1}} \sup_{\nu\in \cT}
    \dist_{\Ptwo}(\nu,\tcB(\mu_0,\dots,\mu_N)) \geq \frac{1}{2\sqrt{3}\,\log(N+1)} \qquad\text{for all }N\in\N .
 \end{equation}
\end{proposition}

\begin{proof}
Let $\cT \coloneqq \{\mu_n : n \in \N_0\} \subset \PtwoR$ 
be the family \eqref{eq:sharp-family-measures} determined by the quantile functions
\begin{equation*}
    \stepfunc_n \coloneqq \Qfunc{\mu_n} = c_n 2^{n/2}\, \bm{1}_{[t_n,1]} 
    \quad \text{with} \quad
    c_n \coloneqq \frac{1}{\log(n+1)}.
\end{equation*}
Since $(c_n)$ is decreasing, Lemma \ref{lem:spike-fam-measure} guarantees that $\cT$ 
is compact, and \eqref{eq:spike-lower-measure} implies
\begin{equation*} 
    \max_{\mu \in \cT}
    \dist_{\Ptwo}\left(\mu,\tcB(\nu_0,\dots,\nu_N)\right)
    \ \geq\ \frac{3^{-1/2}}{\log(N+3)}
    \ \geq\ \frac{3^{-1/2}}{2\log(N+1)}
\end{equation*}
for every $N \ge 1$ and every selection $\nu_0,\dots,\nu_{N} \in \cT$.
In particular,
 $\tsigma_N \geq (2\sqrt{3}\,\log(N+1))^{-1}$
 for every weak greedy sequence for $\cT$.

We now construct six fixed elements  $p_0,\ldots, p_5$, from $a_0,a_1,a_2, q_*\in \Qcone$ defined below, as
\begin{equation*}
    p_0 \coloneqq 0, \qquad p_1 \coloneqq q_*, \qquad p_2 \coloneqq q_*+1, \qquad p_{\ell+3} \coloneqq
    q_*+a_\ell \quad (\ell \in \{0,1,2\}),
\end{equation*}
which, as we show below, belong to $\Qcone$ and their projected affine hull
contains all of $\cK \coloneqq \{\Qfunc{\mu} : \mu \in \cT\}$. 
For $\ell = 0,1,2$, $a_\ell$ is a step function with values denoted by $s_{\ell,j}$ on 
intervals $(t_{j-1},t_j)$ for $j\ge 1$ and defined as
\[
    a_\ell \coloneqq \Phi'_\ell,
\]
with $\Phi_\ell : [0,1] \to \R$ continuous, affine on every
interval $[t_{j-1},t_j]$, and determined by
\begin{equation} \label{eq:Phi-knots}
    \Phi_\ell(0) = \Phi_\ell(1) = 0, \qquad \Phi_\ell(t_j) = j^\ell2^{-j^2}
    \quad (j \ge 1),
\end{equation}
represented in Figure~\ref{fig:phi}.
Since $t_j - t_{j-1} = 2^{-j}$, there exists $c >0$ such that
\[
    |s_{\ell,j}| \le c\, (1+j)^22^{-j^2+3j},
\]
hence both $\sum_j (t_j-t_{j-1})|s_{\ell,j}|^2$ and
$\sum_j|s_{\ell,j}|$ converge, so that $a_\ell \in L^2(0,1)$ with absolutely
summable jumps.

Finally we define $q_*$ as the nondecreasing step function $q_*$
that vanishes on $(t_0,t_1)$ and has jumps at $t_j$
\begin{equation*}
    J_j \coloneqq \sum_{\ell=0}^2|s_{\ell,j+1}-s_{\ell,j}|, \quad j \geq 1.
\end{equation*}
Then, $q_* + a_\ell$ has jump $J_j + s_{\ell,j+1}-s_{\ell,j} \ge 0$ at $t_j$ and 
is therefore nondecreasing, too. Consequently, the six functions $p_0,\ldots, p_5$ 
all belong to $\Qcone$.

We now prove that the projected affine hull of $(p_0,\ldots, p_5)$ contains all of $\cK$.
To do so, we construct functions $h_n \in   \aff\{p_0,\dots,p_5 \}$ below satisfying 
$P_\Qcone (h_n) = \stepfunc_n $.

Fix $n \ge 1$, set $\Dn := c_nt_n\sqrt{1-t_n}$ and $\kappa_n := 4\Dn2^{n^2}$,
and consider
\begin{equation} \label{eq:counter2-Hn}
    H_n := \kappa_n\big( (n^2-\tfrac14)\Phi_0 - 2n\Phi_1 + \Phi_2 \big)
    \ \in\ \linspan\{\Phi_0,\Phi_1,\Phi_2\},
\end{equation}
see Figure~\ref{subfig:HnGn}.
By \eqref{eq:Phi-knots}, its values at the knots are
\begin{equation} \label{eq:counter2-Hn-knots}
    H_n(t_j) = \kappa_n2^{-j^2}\big((n-j)^2-\tfrac14\big), \quad j\ge 1,
\end{equation}
so that $H_n(0) = H_n(1) = 0$, $H_n(t_n) = -\Dn$, and $H_n(t_j) > 0$ for every
$j \ne n$. The quadratic factor $(n-j)^2-\tfrac14$ thus acts as a selector
placing a single negative pit at the prescribed knot $t_n$, without leaving the
fixed three-dimensional space in \eqref{eq:counter2-Hn}.

Next, let $G_n$ be the triangular function joining $(0,0)$, $(t_n,-\Dn)$ and
$(1,0)$. Its left slope is negative and its right slope is positive, so $G_n$ is
convex, and $G_n \le H_n$, because this holds at every knot by
\eqref{eq:counter2-Hn-knots} and both functions are affine in between.
Conversely, every convex $F \le H_n$ satisfies $F(0) \le 0$, $F(t_n) \le -\Dn$
and $F(1) \le 0$, hence lies below the chords joining these values, that is,
below $G_n$. Therefore $G_n$ is the greatest convex minorant of $H_n$.

Finally, define $h_n := H'_n + \Dn/t_n \in \linspan\{1,a_0,a_1,a_2\} \subset \aff\{p_0,\dots,p_5 \}$, 
noting that $p_0 = 0$ (see Figure~\ref{subfig:hn-Phn}). Adding an affine function to a
cumulative graph adds it to the greatest convex minorant, so that
\begin{equation*}
    P_\Qcone (h_n) = \Big(G_n(t)+\frac{\Dn}{t_n}t\Big)'
    = \frac{\Dn}{t_n(1-t_n)}\bm{1}_{[t_n,1]}
    = \frac{c_n}{\sqrt{1-t_n}}\bm{1}_{[t_n,1]}
    = \stepfunc_n .
\end{equation*}
Since also $\stepfunc_0 = 0 = P_\Qcone (p_0)$, we obtain $\cK 
\subset P_\Qcone(\aff\{p_0,\dots,p_5\})$, and hence, $\td_N$ being nonincreasing in $N$,
\begin{equation*}
    \td_N(\cT;\PtwoR) = 0, \qquad N \ge 5 .
\end{equation*}
This concludes the proof.
\end{proof}

In conclusion, $\td_N(\cT;\PtwoR)$ decays faster than any algebraic or exponential
rate, while no bound $\tsigma_N \le C_2(N+1)^{-\beta}$ with $\beta>0$ can hold by
\eqref{eq:sigmaNlowerbound}. The obstruction exploited by this counterexample
is intrinsic to the nonlinearity of $P_\Qcone$: the image of one
fixed finite-dimensional affine space may contain quantile functions with infinitely many
distinct breakpoints, so that a few well-chosen centers $p_i \notin \cK$
reproduce all of $\cK$ exactly. Affine combinations of finitely many
selected elements from $\cK$, in contrast, only ever see the selected
breakpoints, and by Lemma \ref{lma:projection-partition} the projection cannot
create a missing one. It is precisely this nonlinear active-set effect that the
ordinary barycentric width $d_N$ in Theorem \ref{thm:upper} rules out.

\section{Conclusion}
In this article, we systematically investigated the convergence of the greedy algorithm for metric reduced-order modeling using barycenters in $\PtwoR$.
Having already seen in the introduction that ordinary barycenters are insufficient, we were able to achieve the optimal convergence rate of the barycentric $N$-width in the greedy approximation error by employing signed barycenters. 
Subsequently, using another, technically more complex example, we showed that the algebraic relation can be arbitrarily poor without the use of signed barycenters.
Finally, in a last example, we showed that the approximation rates truly coincide only for this combination of classical and signed barycenters.
This example is nevertheless of considerable structural interest, as
it demonstrates that by abandoning a linear vector space and employing the resulting nonlinear projection in the definition of signed barycenters, the connection between greedy error and $N$-width can break.
Interestingly, numerical experiments in~\cite{dalery2026nonlinear} already
indicated that allowing non-convex weights is essential both for accurate
approximation and for extrapolation.

So far, our analysis is restricted to the one-dimensional 2-Wasserstein space, which already raises substantial difficulties.
In particular, the projection onto the cone of quantile functions is inherently nonlinear, which is why our results differ from those obtained in Hilbert and Banach spaces.

We leave for future work the extension of such results to  generic metric spaces, starting with $\Ptwo(\R^d)$ for $d>1$, as in this case, the isometry with the $L^2$-space of quantile functions is lost and very little is known about signed barycenters beyond $d=1$~\cite{Tornabene2025}. 

\section*{Acknowledgements}
The authors thank Marco Veneroni for valuable advice that helped shape the
direction of this work. CM acknowledges the support of the French National Research Agency (ANR) through
grant ANR-24-CPJ1-0099-01, and GD and MP acknowledge the support of ANR through the project NUMERIQ (ANR-24-CE46-2255).
CM and GD acknowledge support from the CNRS through an International Emerging
Action (IEA).

\section*{AI disclosure}

Chatgpt (GPT-5.6 Sol, GPT-6 Astra) and Claude (Opus 5.5) were used for
mathematical exploration and to develop some of the proofs during this project. The authors developed, completed and verified all arguments and take full responsibility for the results.

\bibliographystyle{amsalpha}
\bibliography{biblio.bib}

\appendix
\section{
  Upper and lower bounds on the widths of families of 
  sums of Dirac masses
}
\label{sec:step-fun}
In this section, we prove the properties of the family \eqref{eq:sharp-family-measures} from Lemma \ref{lem:spike-fam-measure} by making the connection to the quantile functions.
\begin{lemma}
Recalling \eqref{eq:sharp-family}, for $n \in \N$ consider the sequence of breakpoints $t_n = 1-2^{-n}$, and a
sequence of positive coefficients $c_n > 0$ with $c_n \to 0$ to generate the family
\begin{equation*} 
    \stepfunc_n = c_n\, 2^{n/2}\,\bm{1}_{[t_n,1]},
\end{equation*}
as well as $\stepfunc_0 = 0 \in L^2(0,1)$.
Then $\stepfunc_n \in \Qcone$ with $\|\stepfunc_n\| = c_n$, and the set $\cK \coloneqq \{\stepfunc_n : n \in \N_0\}$
is compact. 
Moreover, for all $N\in \N_0$,
\begin{equation}
  \label{eq:N-width-step}
  d_N(\cK; \Qcone) \,\leq\, \max_{n\geq N+1} c_n,
\end{equation}
and for every $N \ge 0$ and every selection $f_0, \dots , f_N \in
\cK$,
\begin{equation} \label{eq:spike-lower}
    \sup_{n \ge 0}\ \dist_{L^2}\Big(\stepfunc_n,
    P_\Qcone (\aff\{f_0,\dots,f_{N}\})\Big)
    \, \geq\, \min_{1\le n\le N+2} \frac{c_n}{\sqrt{3}} .
\end{equation}
\end{lemma}
\begin{proof}
Each $\stepfunc_n$ is nondecreasing and bounded, hence $\stepfunc_n \in \Qcone$, and from
$(1-t_n) = 2^{-n}$ it follows that
\[
    \|\stepfunc_n\|^2 = c_n^2 \longrightarrow 0,
\]
so that $\cK$ is compact.

For $0 \le n \le N$ we have $\stepfunc_n \in \cB(\stepfunc_0,\dots,\stepfunc_N)$, while
for $n \geq N+1$, using $\stepfunc_0 \in \cB(\stepfunc_0,\dots,\stepfunc_N)$,
\begin{equation*}
    \dist_{L^2}\bigl(\stepfunc_n,\cB(\stepfunc_0,\dots,\stepfunc_N)\bigr) \le \|\stepfunc_n - \stepfunc_0\| =
    c_n,
\end{equation*}
which gives \eqref{eq:N-width-step}.

Then, fix $N \ge 0$ and a selection $f_0,\dots,f_{N} \in \cK$, and let $S \subset \N$ denote the set of positive indices occurring in it 
(i.e., the indices $i\in\N$ such that there exists  $\stepfunc_i\in \{f_0, \dots, f_{N}\}$), 
so that $\# S \le N+1$ and
\[
   A := \aff\{f_0,\dots,f_{N}\} \subset\aff\big(\{\stepfunc_j : j \in S\} \cup
   \{\stepfunc_0\}\big).
\]
Let $R$ denote the partition of $(0,1)$ whose interior breakpoints are $\{t_j
: j \in S\}$, and let $V_R$ be as in Lemma \ref{lma:projection-partition}.
Every $\stepfunc_j$ with $j \in S$ is constant on each subinterval of $R$, and so is
$\stepfunc_0 = 0$; hence $A \subset V_R$, and Lemma \ref{lma:projection-partition}
gives
\begin{equation} \label{eq:spike-in-VR}
    P_\Qcone(A) \subset V_R .
\end{equation}

Since $\# S \le N+1$, at least one index $n \in \{1,\dots,N+2\}$ does not belong
to $S$.
The missing breakpoint $t_n$ then lies strictly inside a cell $(u,v)$ of $R$.
As $t_i < t_j$ whenever $i < j$, the neighboring breakpoints are at worst
$t_{n-1}$ and $t_{n+1}$, so that
\[
    \ell_1 := t_n - u \geq 2^{-n},
    \qquad
    \ell_2 := v-t_n \geq 2^{-n-1}.
\]
Any $h \in V_R$ is constant on $(u,v)$, say equal to $c$, while $\stepfunc_n = 0$ on
$(u,t_n)$ and $\stepfunc_n = a := c_n2^{n/2}$ on $[t_n,v)$. It follows
\begin{equation}
    \|\stepfunc_n - h\|^2
    \ \ge\ \min_{c \in \R}\big( \ell_1c^2 + \ell_2(a-c)^2 \big)
    = a^2\frac{\ell_1\ell_2}{\ell_1+\ell_2}
    \ \ge\ a^2\frac{2^{-n}}{3}
    = \frac{c_n^2}{3},
\end{equation}
and \eqref{eq:spike-lower} follows together with \eqref{eq:spike-in-VR}.
\end{proof}

\begin{proof}[Proof of Lemma~\ref{lem:spike-fam-measure}]
For $\mu_0 = \delta_0$ we have $\Qfunc{\mu_0} = \stepfunc_0 = 0$.
For $n \ge 1$ 
the cumulative distribution function of $\mu_n$ is 
\[ 
  \mu_n((-\infty,x]) = \begin{cases}
    0 & x \in (-\infty,0) \\ 
    t_n & x \in [0,c_n 2^{n/2}) \\ 
    1 & x \in [c_n 2^{n/2}, \infty)
  \end{cases} 
\]
Hence $\Qfunc{\mu_n}(t) = 0$ for $t \in (0,t_n]$ and $\Qfunc{\mu_n}(t) = c_n 2^{n/2}$ for $t \in (t_n,1)$, 
i.e.~$\Qfunc{\mu_n} = \stepfunc_n$ in $L^2(0,1)$.
Since $\mu \mapsto \Qfunc{\mu}$ is a bijective isometry from $(\PtwoR,\dist_{\Ptwo})$ 
onto $(\Qcone,\dist_{L^2})$, it maps $\cT \subset \PtwoR$ isometrically onto $\cK \subset L^2(0,1)$.

By Corollary~\ref{cor:bary-quant-character}, for every
$\bm\lambda \in \tLambda_N$ the signed barycenter is unique and for all $\nu_0,\dots,\nu_N \in \PtwoR$
\[ 
  \Qfunc{\bary{\bm\lambda}{\nu_0,\dots,\nu_N}} 
  = P_\Qcone\left(\textstyle \sum\limits_{i=0}^N \lambda_i \Qfunc{\nu_i}\right).
\]
We now prove the three claims of Lemma~\ref{lem:spike-fam-measure}.

(i) $\cT$ is the preimage of the compact set $\cK$ under the isometric bijection 
$\mu \mapsto \Qfunc{\mu}$, so it is compact.

(ii) Let $N \geq 0$ and let $\nu_0,\dots,\nu_N \in \PtwoR$. Letting $\bm\lambda$ 
range over $\Lambda_N$ obtain
\[ 
  \Qfunc{\cB(\nu_0,\dots,\nu_N)} = \conv\{\Qfunc{\nu_0},\dots,\Qfunc{\nu_N}\}. 
\]
Since $\mu \mapsto \Qfunc{\mu}$ is onto $\Qcone$, it follows that 
$d_N(\cT;\PtwoR) = d_N(\cK;\Qcone)$, and \eqref{eq:spike-width-measure} is 
\eqref{eq:N-width-step}.

(iii) Let $N \ge 0$, let $\nu_0,\dots,\nu_{N} \in \cT$, and set
$f_i := \Qfunc{\nu_i} \in \cK$. 
Letting $\bm\lambda$ range over $\tLambda_{N}$ obtain
\[
  \Qfunc{\tcB(\nu_0,\dots,\nu_{N})} = P_\Qcone\bigl(
    \aff\{\Qfunc{\nu_0},\dots,\Qfunc{\nu_{N}}\}
  \bigr),
\]
and for every $n \ge 0$,
\[
  \dist_{\Ptwo}\bigl(\mu_n,\tcB(\nu_0,\dots,\nu_{N})\bigr)
  = \dist_{L^2}\bigl(\stepfunc_n, P_\Qcone(\aff\{f_0,\dots,f_{N}\})\bigr).
\]
Taking the supremum over $n$ and applying \eqref{eq:spike-lower} gives
\eqref{eq:spike-lower-measure}.
\end{proof} 

\end{document}

%% file: tikz-example1.tex
\begin{tikzpicture}
\pgfmathsetmacro{\cx}{1.35}    
\pgfmathsetmacro{\cy}{0.43}   
\pgfmathsetmacro{\rad}{1}     
\pgfmathsetmacro{\rot}{0}     
\pgfmathsetmacro{\side}{2*sqrt(3)*\rad+0.25} 

\pgfmathsetmacro{\Rtri}{\side/sqrt(3)} 
\pgfmathsetmacro{\xmin}{min(-0.2, \cx-\Rtri/2-0.2)}
\pgfmathsetmacro{\xmax}{\cx+\Rtri+0.5}
\pgfmathsetmacro{\ymin}{min(-0.2, \cy-\side/2-0.2)}
\pgfmathsetmacro{\ymax}{max( 0.2, \cy+\side/2+0.21}

\definecolor{wong1}{HTML}{0072b2}
\definecolor{wong2}{HTML}{e69f00}
\definecolor{wong3}{HTML}{009e73}
\begin{axis}[
	x=1.8cm,y=1.8cm,
	axis lines=middle,
    xlabel style={anchor=north east,},
    ylabel style={anchor=north east,},
	xmin=\xmin, xmax=\xmax,
	ymin=\ymin, ymax=\ymax,
	xtick distance=0.5, ytick distance=0.5,
	xticklabel=\empty, yticklabel=\empty,
	xlabel=$s$, ylabel=$m$,
	line cap=round
]

\coordinate (C) at (\cx,\cy);
\pgfplotsinvokeforeach{0,...,4}{
	\coordinate (Q#1) at ({\cx+\rad*cos(\rot+72*#1)}, {\cy+\rad*sin(\rot+72*#1)});
}
\coordinate (T1) at ({\cx-\Rtri/2}, {\cy+\side/2});
\coordinate (T2) at ({\cx-\Rtri/2}, {\cy-\side/2});
\coordinate (T3) at ({\cx+\Rtri},   {\cy});
\coordinate (S) at ({\cx+\rad*cos(\rot+36)}, {\cy+\rad*sin(\rot+36)});
\coordinate (M) at ($(Q0)!0.5!(Q1)$);
\coordinate (L) at ({\cx+0.89*\rad*cos(\rot+122)}, {\cy+0.89*\rad*sin(\rot+122)});

\fill[color=wong3, fill opacity=0.1] (T1) -- (T2) -- (T3) -- cycle;

\draw [line width=1.3pt, color=wong2, fill=wong2, fill opacity=0.2] (C) circle (\rad);

\fill[color=wong1, fill opacity=0.8] (Q0) -- (Q1) -- (Q2) -- (Q3) -- (Q4) -- cycle;

\draw [line width=1pt, color=wong1] (Q0) -- (Q1) -- (Q2) -- (Q3) -- (Q4) -- cycle;
\draw [line width=1pt, color=wong3] (T1) -- (T2) -- (T3) -- cycle;
\node[wong3, above] at (T1) {$\nu_0$};
\node[wong3, below] at (T2) {$\nu_1$};
\node[wong3, right] at (T3) {$\nu_2$};

\draw [{To[length=2pt]}-{To[length=2pt]}, line width=0.8pt, color=black] (S) -- (M);

\node[pin=96:{\color{wong2}{$\mathcal{T}$}}, color=black] at (L) {};
\path (T2) -- (T3) node[pos=0.6, below right, color=wong3] {$d_2=0$};
\node[color=white] at (C) {$\cB(\mu_0,\dots,\mu_4)$};

\node[pin=85:{{\color{black}
   $\sigma_4 > 0$ 
   }}] at (S) {};

\end{axis}
\end{tikzpicture}

%% file: phi-functions.tex
\begin{tikzpicture}[
  x=10cm,y=10cm,
  font=\small,
  line cap=round,line join=round
]
  \definecolor{figOrange}{HTML}{B44924}
  \definecolor{figBlue}{HTML}{246595}
  \definecolor{figGreen}{HTML}{347C56}
  \definecolor{figInk}{HTML}{242424}

  \foreach \v/\lab in {
    0/0,0.125/{\frac18},0.25/{\frac14},
    0.375/{\frac38},0.5/{\frac12}
  }{
    \ifdim\v pt>0pt
      \draw[gray!18,line width=0.4pt] (0,\v) -- (1,\v);
    \fi
    \draw[figInk,line width=0.4pt] (-0.006,\v) -- (0,\v);
    \node[anchor=east,inner sep=4pt] at (0,\v) {$\lab$};
  }
  \draw[figInk,line width=0.5pt] (0,0.55) -- (0,0) -- (1.025,0);
  \node[anchor=west,inner sep=3pt] at (1.025,0) {$t$};
  \foreach \x/\lab in {0/0,0.5/t_1,0.75/t_2,0.875/t_3,1/1}{
    \draw[figInk,line width=0.4pt] (\x,0) -- (\x,-0.007);
    \node[anchor=north,inner sep=4pt] at (\x,0) {$\lab$};
  }

  \draw[figInk,line width=1.05pt] (0,0) -- (0.5,0.5);
  \draw[figBlue,line width=1.05pt] plot coordinates {
    (0.5,0.5) (0.75,0.0625) (0.875,0.001953125)
    (0.9375,0.0000152587890625)
    (0.96875,0.0000000298023224)
    (0.984375,0.0000000000145519) (1,0)
  };
  \draw[figOrange,line width=1.05pt] plot coordinates {
    (0.5,0.5) (0.75,0.125) (0.875,0.005859375)
    (0.9375,0.00006103515625)
    (0.96875,0.0000001490116119)
    (0.984375,0.0000000000873115) (1,0)
  };
  \draw[figGreen,line width=1.05pt] plot coordinates {
    (0.5,0.5) (0.75,0.25) (0.875,0.017578125)
    (0.9375,0.000244140625)
    (0.96875,0.0000007450580597)
    (0.984375,0.0000000005238689) (1,0)
  };

  \foreach \x/\y in {
    0.75/0.0625,0.875/0.001953125,0.9375/0.0000152587890625
  }{\fill[figBlue] (\x,\y) circle[radius=1.2pt];}
  \foreach \x/\y in {
    0.75/0.125,0.875/0.005859375,0.9375/0.00006103515625
  }{\fill[figOrange] (\x,\y) circle[radius=1.2pt];}
  \foreach \x/\y in {
    0.75/0.25,0.875/0.017578125,0.9375/0.000244140625
  }{\fill[figGreen] (\x,\y) circle[radius=1.2pt];}
  \foreach \x/\y in {0/0,0.5/0.5,1/0}{
    \fill[figInk] (\x,\y) circle[radius=1.5pt];
  }

  \foreach \y/\col/\ell in {
    0.50/figBlue/0,0.44/figOrange/1,0.38/figGreen/2
  }{
    \draw[\col,line width=1.1pt] (0.88,\y) -- (0.95,\y);
    \node[\col,anchor=west,inner sep=4pt] at (0.95,\y) {$\Phi_\ell$};
  }
\end{tikzpicture}

%% file: hn-gn.tex
\begin{tikzpicture}[
  x=12cm,y=1.55cm,
  font=\small,
  line cap=round,line join=round,
  >=Stealth
]
  \definecolor{figOrange}{HTML}{B44924}
  \definecolor{figBlue}{HTML}{246595}
  \definecolor{figInk}{HTML}{242424}
  \definecolor{figFill}{HTML}{EDF2F6}

  \def\Hcoordinates{
    (0,0)
    (0.5,1.65)
    (0.75,0.82)
    (0.875,-1)
    (0.9375,0.27)
    (0.96875,0.10)
    (0.984375,0.034)
    (0.9921875,0.011)
    (0.99609375,0.0035)
    (0.998046875,0.001)
    (1,0)
  }

  \fill[figFill]
    plot coordinates {\Hcoordinates}
    -- (0.875,-1) -- (0,0) -- cycle;

  \draw[gray!65,densely dashed,line width=0.4pt]
    (0,-1) -- (0.875,-1);
  \draw[gray!65,densely dashed,line width=0.4pt]
    (0.875,0) -- (0.875,-1);
  \draw[figInk,->,line width=0.5pt] (-0.025,0) -- (1.045,0);
  \draw[figInk,->,line width=0.5pt] (0,-1.30) -- (0,1.97);

  \draw[figBlue,line width=1.05pt]
    (0,0) -- (0.875,-1) -- (1,0);
  \draw[figOrange,line width=1.1pt]
    plot coordinates {\Hcoordinates};

  \foreach \x/\y in {
    0.5/1.65,0.75/0.82,0.9375/0.27,0.96875/0.10
  }{
    \fill[figOrange] (\x,\y) circle[radius=1.2pt];
  }
  \foreach \x/\y in {0/0,0.875/-1,1/0}{
    \fill[figInk] (\x,\y) circle[radius=1.5pt];
  }
  \foreach \x in {0.5,0.75,0.875,0.9375,1}{
    \draw[figInk,line width=0.4pt] (\x,-0.025) -- (\x,0.025);
  }

  \node[anchor=north east,inner sep=2pt] at (0,0) {$0$};
  \node[anchor=north,inner sep=2pt] at (1,0) {$1$};
  \node[anchor=west,inner sep=2pt] at (1.045,0) {$t$};
  \node[anchor=east,inner sep=3pt] at (0,-1) {$-\Dn$};
  \node[anchor=north,inner sep=3pt] at (0.5,0) {$t_1$};
  \node[anchor=north,inner sep=3pt] at (0.75,0) {$t_{n-1}$};
  \node[anchor=south,inner sep=3pt] at (0.875,0) {$t_n$};
  
  \node[figOrange,font=\normalsize,anchor=south] at (0.59,1.50) {$H_n$};
  \node[figBlue,font=\normalsize] at (0.39,-0.72) {$G_n$};
  \node[anchor=north,inner sep=4pt] at (0.875,-1.08) {$(t_n,-\Dn)$};
\end{tikzpicture}

%% file: hn-projection.tex
\begin{tikzpicture}[
  x=12cm,y=0.9cm,
  font=\small,
  line cap=round,line join=round,
  >=Stealth
]
  \definecolor{figOrange}{HTML}{B44924}
  \definecolor{figBlue}{HTML}{246595}
  \definecolor{figInk}{HTML}{242424}
  \definecolor{figFill}{HTML}{EDF2F6}

  \fill[figFill] (0.875,0) rectangle (1,1);
  \draw[figInk,->,line width=0.5pt] (-0.025,0) -- (1.045,0);
  \draw[figInk,->,line width=0.5pt] (0,-3.85) -- (0,3.70);

  \foreach \x/\yBefore/\yAfter in {
    0.5/3.1804912607/-3.3109628413,
    0.75/-3.3109628413/-2.73219937,
    0.875/-2.73219937/1.916015625,
    0.9375/1.916015625/0.0437698364,
    0.96875/0.0437698364/0.1233996525,
    0.984375/0.1233996525/0.12499635,
    0.9921875/0.12499635/0.1249999984,
    0.99609375/0.1249999984/0.125,
    0.998046875/0.125/0.125,
    0.9990234375/0.125/0.125,
    0.9995117188/0.125/0.125,
    0.9997558594/0.125/0.125,
    0.9998779297/0.125/0.125,
    0.9999389648/0.125/0.125
  }{
    \draw[figOrange,opacity=0.60,densely dashed,line width=0.5pt]
      (\x,\yBefore) -- (\x,\yAfter);
  }
  \foreach \a/\b/\z in {
    0/0.5/3.1804912607,
    0.5/0.75/-3.3109628413,
    0.75/0.875/-2.73219937,
    0.875/0.9375/1.916015625,
    0.9375/0.96875/0.0437698364,
    0.96875/0.984375/0.1233996525,
    0.984375/0.9921875/0.12499635,
    0.9921875/0.99609375/0.1249999984,
    0.99609375/0.998046875/0.125,
    0.998046875/0.9990234375/0.125,
    0.9990234375/0.9995117188/0.125,
    0.9995117188/0.9997558594/0.125,
    0.9997558594/0.9998779297/0.125,
    0.9998779297/0.9999389648/0.125,
    0.9999389648/1/0.125
  }{
    \draw[figOrange,line width=1.1pt] (\a,\z) -- (\b,\z);
  }

  \draw[figBlue,line width=1.2pt] (0,0) -- (0.875,0);
  \draw[figBlue,densely dashed,line width=0.7pt] (0.875,0) -- (0.875,1);
  \draw[figBlue,line width=1.2pt] (0.875,1) -- (1,1);

  \foreach \x/\lab in {
    0/0,0.5/t_1,0.75/{t_{n-1}},0.875/t_n,0.9375/{t_{n+1}},1/1
  }{
    \draw[figInk,line width=0.4pt] (\x,-0.045) -- (\x,0.045);
    \node[anchor=north east,inner sep=2pt] at (\x,0) {$\lab$};
  }
  \node[anchor=west,inner sep=3pt] at (1.045,0) {$t$};
  
  \draw[figOrange,line width=1.1pt] (0.16,-2.68) -- (0.215,-2.68);
  \node[figOrange,anchor=west,inner sep=4pt,font=\normalsize] at (0.215,-2.68) {$h_n$};
  \draw[figBlue,line width=1.2pt] (0.16,-3.48) -- (0.215,-3.48);
  \node[figBlue,anchor=west,inner sep=4pt,font=\normalsize] at (0.215,-3.48) {$P_{\Qcone}h_n$};
\end{tikzpicture}